%% file: ms.tex
\documentclass[10pt]{article} 

\usepackage{epstopdf}
\usepackage{algorithmic}
\usepackage{tabularx}
\usepackage{longtable}
\usepackage{graphicx}

\usepackage{amssymb,amsmath,amsfonts,amsthm}
\usepackage{thmtools}

\usepackage{alltt}
\usepackage{xcolor}
\usepackage{float}
\usepackage{colortbl}
\usepackage{arydshln}
\usepackage[font=small,labelfont=bf]{caption}
\usepackage{mathtools}
\usepackage{diagbox}
\usepackage[shortlabels]{enumitem}
\usepackage{bbold}
\usepackage[T1]{fontenc}
\usepackage{hyperref}
\usepackage{nameref,cleveref}
\usepackage{tikz}
\usepackage{csquotes}
\usepackage{nicefrac}

\usepackage{todonotes}
\usepackage{geometry}
\usepackage{mathrsfs}
 
\usepackage{mathrsfs}  
\usepackage{enumitem}
\usepackage{pifont}

\title{Sharp-interface limits of Cahn--Hilliard models and mechanics with moving contact lines}

\author{Leonie Schmeller\thanks{Weierstrass Institute, Mohrenstrasse 39, 10117 Berlin, Germany 
  (\href{mailto:leonie.schmeller@wias-berlin.de}{leonie.schmeller@wias-berlin.de}, \href{mailto:dirk.peschka@wias-berlin.de}{dirk.peschka@wias-berlin.de})}
\and Dirk Peschka\footnotemark[1]}

\usepackage{amsopn}

\input{marcos.tex}

\begin{document}
\maketitle
\begin{abstract}
We construct gradient structures for free boundary problems with nonlinear elasticity and study the impact of moving contact lines. In this context, we numerically analyze how phase-field models converge to certain sharp-interface limits when the interface thickness tends to zero $\varepsilon\to 0$. In particular, we study the scaling of the Cahn-Hilliard mobility $m(\varepsilon)=m_0\varepsilon^\alpha$ for $0\le \alpha \le \infty$. In the presence of interfaces, it is known that the intended sharp-interface limit is only valid for $\underline{\alpha}<\alpha<\overline{\alpha}$.  However, in the presence of moving contact lines we show that  $\alpha$ near $\underline{\alpha}$ produces significant errors.
\end{abstract}

\input{input_text.tex}

\section*{Acknowledgement}
Both authors thank Marita Thomas, for discussions about analysis and nonlinear
elasticity and Andreas M\"unch regarding aspects of asymptotics of phase-field models with degenerate mobilities. 
Furthermore, we thank Helmut Abels and
Harald Garcke for the discussion on sharp-interface limits and lower and upper bounds
of the mobility.

Both authors acknowledge the funding by the German Research Foundation (DFG) within the DFG Priority Program SPP 2171 \emph{Dynamic Wetting of Flexible, Adaptive, and Switchable Substrates} through the projects \#422786086 (LS) and \#422792530 (DP). LS \& DP thank the Berlin Mathematics Research Center MATH$^+$ for funding and support through project AA2-9.

\bibliographystyle{abbrv}
\bibliography{ms}
\end{document}

%% file: marcos.tex
\newcommand{\bu}{\boldsymbol{u}}
\newcommand{\bchi}{\boldsymbol{\chi}}

\newcommand{\bF}{\boldsymbol{F}}

\newcommand{\Wel}{W_{\rm elast}}

\newcommand{\Wpfe}{W_{\rm phase}^{\varepsilon}}
\newcommand{\mob}{m}

\newcommand{\R}{\mathbb{R}}

\newcommand{\dx}{\,{\rm d}x}

\newcommand{\calF}{\mathscr{F}}

\newcommand{\qd}{q_{\varepsilon}}
\newcommand{\qs}{q_0}

\newcommand{\calVd}{\mathcal{V}_{\varepsilon}}

\newcommand{\EEE}{\color{black}}

\DeclareMathOperator*{\esssup}{ess\,sup}
\newtheorem{remark}{Remark}

\crefalias{subsection}{section}

%% file: input_text.tex

\section{Introduction}
Interfaces and surfaces with surface energy are ubiquitous in nature, and their description can employ different levels of detail. 
{For multiphase systems with moving interfaces different mathematical and numerical techniques are available. Explicit descriptions via front tracking methods \cite{tryggvason2001front} allow high control over the moving boundary.
Implicit description via level-set methods \cite{chang1996level,osher1988fronts} or using phase fields lead to diffuse representations of interfaces \cite{anderson1998diffuse}.
The choice of method depends on problem-specific features and requirements, e.g. feasibility of topological changes, treatment of discontinuities at the interface, mass conservation and transport across interfaces, ease of  implementation, availability of computational resources, precision and control over interface evolution.}
In general for continuum models, sharp-interface and phase-field models are the most common levels of abstraction with more or less details, respectively. 

Such systems can be arbitrarily complicated when interfaces have a complex substructure, e.g., with interfacial mass and thermodynamical effects. 
This complexity increases even further by considering nonlinear diffusion and phase separation, reactions and phase transitions and complex viscoelastic properties as possible dissipative processes. 
However, even for simple systems with constant surface tension and fluid flow or elasticity in the pure phase, the connection between common diffuse-interface models and their sharp limits is elusive.
In phase-field models, the interface is modeled by continuous functions $\psi$ that are constant inside a phase, e.g. $\psi=\pm 1$, and whose change in an interfacial region of width proportional to $\varepsilon>0$ indicates the transition to another phase in a continuous manner. The mathematical purpose of the phase field is twofold. Firstly, it should act as an indicator for the presence of a phase. Secondly, it defines a phase-field energy density $\smash{\Wpfe(\psi,\bF^{-T}\nabla\psi)}$ that measures, among other things, the surface tension of the interface with deformation gradient $\bF$. Without mechanics $\bF=\boldsymbol{I}$ the Cahn-Hilliard model 
\begin{align}
\label{eqn:base_cahn_hilliard}
\begin{split}
\partial_t\psi = \nabla\cdot (m \nabla \mu_\varepsilon)\,, \quad \mu_\varepsilon = 
\frac{\delta}{\delta\psi}\Wpfe\,,\quad
\Wpfe=\frac{3}{2\sqrt{2}}\Big[\tfrac\varepsilon2|\nabla\psi|^2+\tfrac1{4\varepsilon}(1-\psi^2)^2\Big],
\end{split}
\end{align}
is commonly used to model phase-field evolution and phase separation. Its well-known that $\smash{\Wpfe}$ converges to the perimeter of the interface connecting regions with $\psi=\pm 1$.
The limit passage of dynamic phase-field models to sharp-interface models can be analyzed by different techniques, e.g. by matched asymptotics \cite{meca2018sharp,lee2016sharp, abels2012thermodynamically}, $\Gamma$-convergence \cite{modica1987gradient} and evolutionary convergence. 
{
Degenerate $\psi$-dependent mobilities $m(\varepsilon,\psi)$ are highly relevant  
for the limit passage from the diffuse to the sharp-interface model, which is studied in detail in \cite{abels2012thermodynamically,meca2018sharp,dziwnik2017anisotropic}.
}

In many studies, only one set of interfacial width $\varepsilon$ and mobility $\mob$ is used for a specific application leading to reasonable results \cite{van2021adaptive, di2014cahn, liu2003phase, boyer2010cahn, mokbel2018phase}. 
Instead, we focus on a systematic study of the appropriate choice of the Cahn-Hilliard mobility. {We extend the work by Yue et al.~\cite{yue2010sharp} by considering fluid-structure interaction and by directly comparing sharp and diffuse-interface models.}
The goal of this work is to establish a general numerical approach to the problem of sharp-interface limits with moving contact lines, somewhat in the spirit of previous work by Liu et al.~\cite{liu2022diffuse} and Aland et al.~\cite{aland2010two}.
{We consider $\mob(\varepsilon)=m_0\varepsilon^{\alpha}$ with an appropriate $0\leq\alpha\leq\infty$ and determine the scaling of $\mob$. 
Usually, for interface problems \footnote{excluding cases with three-phase contact lines}, $m=\varepsilon m_0$ yields the intended interface condition \cite{abels2012thermodynamically} and gives an upper bound on the mobility. 
Further, the authors in \cite{abels2022non, schaubeck2014sharp} show that there is also a lower bound on the Cahn-Hilliard mobility which suggests that the intended limit is reached for $\alpha\in(\underline{\alpha},\bar{\alpha})$.}
{However, since such techniques are inevitably much more complex for moving contact lines and require considerations in higher-dimensions, we will resort to numerical techniques.} 

In order to explore the sharp-interface limit of phase-field models, the strategy of the paper is as follows.
In \Cref{sec:models} the phase-field evolution \eqref{eqn:base_cahn_hilliard} is extended to a ternary multiphase systems with solid (s), liquid ($\ell$) and air (a) phases. Similarly, we provide a sharp-interface model of the supposedly same system.
Nonlinear elasticity of compressible materials and appropriate coupling terms for phase-field and sharp-interface model are introduced using an hyperelastic energy density $W_\mathrm{elast}(\bF,\psi)$ with deformation gradient $\bF=\nabla\bchi$. At the moving contact line we assume an equilibrium of surface forces. 
We  provide an Eulerian formulation of the system of nonlinearly coupled PDEs and also provide the corresponding Lagrangian weak formulation that uses the underlying gradient flow structure of the model.

For the phase-field model we discuss different possible scaling limits of the mobility $m$ in order to have convergence to the sharp-interface model as $\varepsilon\to 0$. We also review existing results on scaling laws and options for choosing the Cahn-Hilliard mobility for models with interfaces. Then we introduce the technical methodology with which we perform our numerical convergence analysis.
In \Cref{sec:si_vs_di} we perform the discretization of both phase-field and sharp-interface model in order to characterize the sharp-interface limit with contact lines using suitable norms $\|\cdot\|$. The goal here is to discuss the distance $\|\bchi_{0,\Delta}-\bchi_{\varepsilon,\Delta'}\|$ of numerical solutions of the sharp-interface model $\bchi_{0,\Delta}$ and of the phase-field model $\bchi_{\varepsilon,\Delta'}$
after interpolating between different discrete function spaces for different discretization parameters $\Delta,\Delta'$. 

Most importantly, this requires a discussion of discretization errors of the two models with respect to time-discretization, mesh refinement and the polynomial order of the finite element spaces. Therefore, as the main tool we use the Bochner norms for  $L^2(Q_T;\mathbb{R}^d)$ and $L^\infty(Q_T;\mathbb{R}^d)$ for the deformation to assess the error in the space-time cylinder $Q_T=[0,T]\times\Omega$ as $\varepsilon\to 0$ for different $m(\varepsilon)=m_0\varepsilon^\alpha$.
Finally, in \Cref{sec:convergence_discussion} we perform a discussion of the numerical solutions generated with the methodology and based on the error analysis introduced in \Cref{sec:si_vs_di}. We  focus on comparing concrete solutions of the benchmark problem ``liquid droplet on viscoelastic substrate''.

The main goal of \Cref{sec:convergence_discussion} and this entire work is to identify exponents $\underline{\alpha}\le\alpha_\mathrm{opt}\le\overline{\alpha}$ that give an optimal mobility $m_\mathrm{opt}=m_0\varepsilon^{\alpha_\mathrm{opt}}$ in the sense of (error of) the sharp-interface limit and to provide a systematic approach to identify such $\alpha_\mathrm{opt}$. While such a statement usually depends on the norm, our choice of norm and discussion is based on the prediction of the moving contact line.

\section{Models for diffuse and sharp interfaces in variational form}
\label{sec:models}

Cahn-Hilliard-Navier-Stokes systems are widely used for modeling multiphase flows in various applications as well as for incompressible viscous Newtonian fluids \cite{boyer2010cahn,abels2012thermodynamically}.
Since the numerical treatment of physical interface thicknesses $\varepsilon$ on the nanoscale is often infeasible for phase-field models, it is important to clearly understand the behavior for $\varepsilon\to 0$, i.e. the \emph{sharp-interface limit}. 
The Cahn-Hilliard mobility $\mob$ depends on $\varepsilon$ and possibly on the phase field variable $\psi$ \cite{novick2008cahn}. An advantage of $\psi$-dependent degenerate mobilities is the ability to guarantee that the phase field variables remain in a fixed, application-specific interval, which in general is a absent property. 

{As a typical test case we use an outer cylinder $\Omega=[0,R]\times [0,H]$ with $R=1$ and $H=1$. The region $\Omega_s=\{(r,z)\in\Omega:0<z<1\}$ defines the solid phase, $\Omega_\ell=\{(r,z)\in\Omega:\sqrt{r^2-(z-1)^2}<r_\text{drop}, z>1\}$ with $r_\text{drop}=\nicefrac12$ defines the liquid phase, and the remainder $\Omega_a=\Omega\setminus(\Omega_\ell\cup\Omega_s)$ is the air phase. The undeformed reference domains $\Omega_s$ (red) and $\Omega_\ell$ (blue) are shown in the left image of \Cref{fig:3d}. The corresponding deformed domains $\bar\Omega_i(t)=\bchi(t,\Omega_i)$ from the sharp-interface model at time $t=T$ are shown in the right image of \Cref{fig:3d}. In this image, the air phase is not shown to improve visibility of the liquid and solid phases and their three-dimensional impression.

\begin{figure}[hb!]\label{fig:3d_drops}
\centering

\includegraphics[width=0.35\textwidth,trim=5cm 6cm 5cm 6cm,clip]{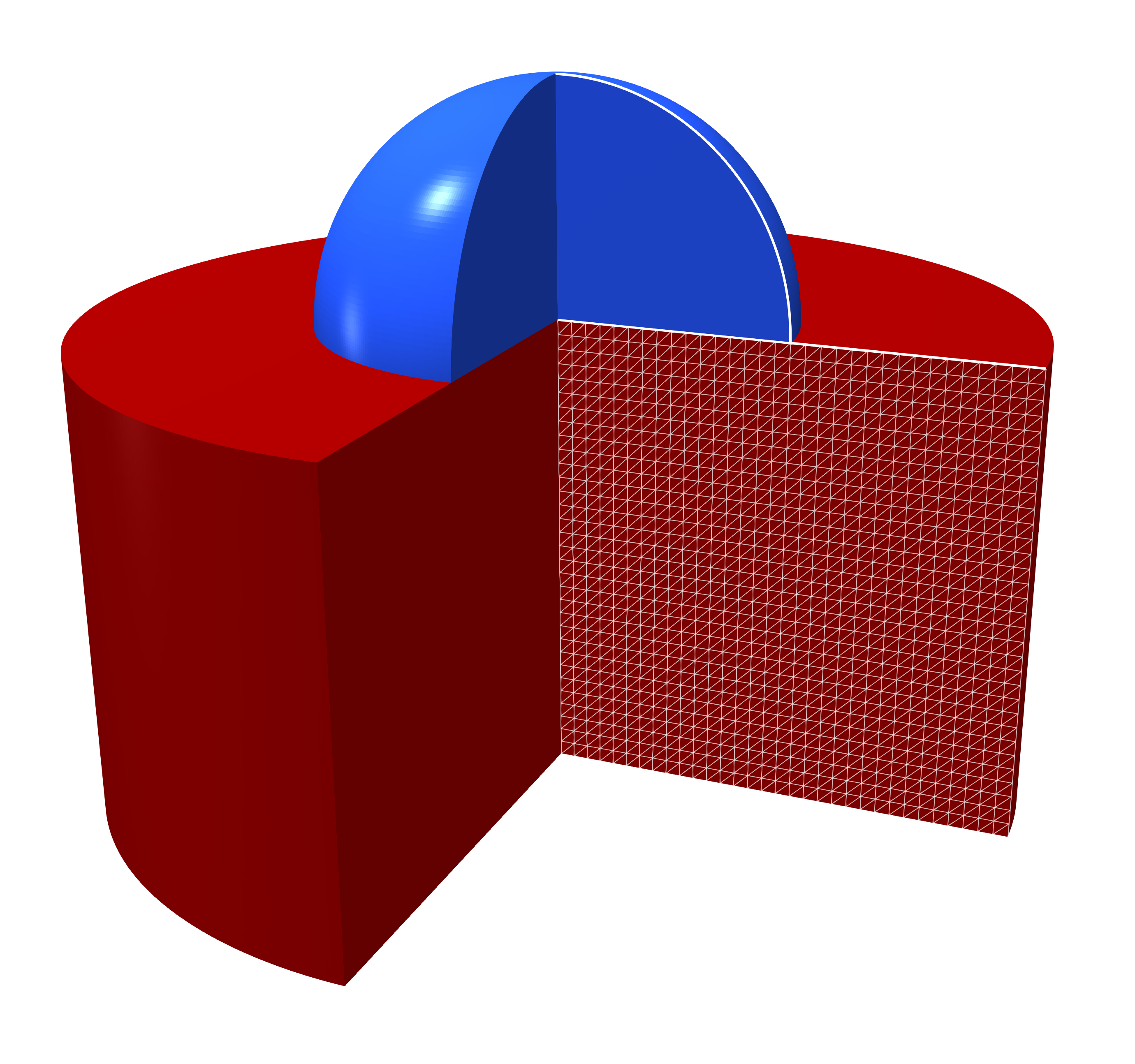}\hspace{1.5cm}
\includegraphics[width=0.35\textwidth,trim=5cm 6cm 5cm 6cm,clip]{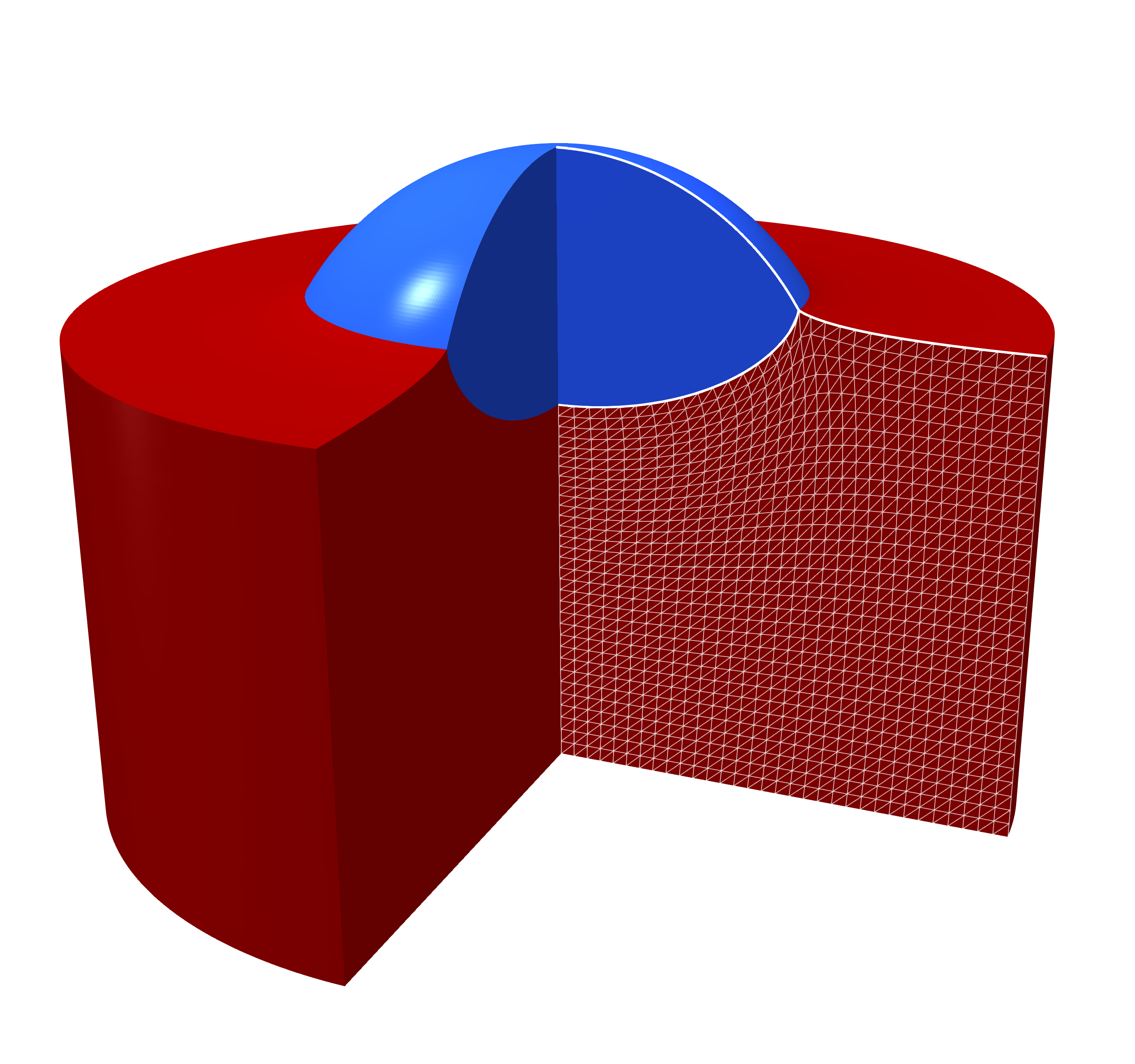}

\caption{Slice of (left) initial liquid droplet (blue half sphere) on a flat elastic substrate (red cylinder). For $x=(r,z)\in\Omega$, at $z=0$ (solid bottom) and $z=H$ (air top) we have no-slip conditions $\bchi(x)=(\chi_r(x),\chi_z(x))=(r,z)$ and at $r=0$ and $r=R$ we fix the normal component $\chi_r(r)=r$ and have natural boundary conditions (symmetry/slip) for the tangential component $\chi_z$. (Right) Deformed droplet $\bar\Omega_\ell(t)$ (blue) and deformed solid substrate $\bar\Omega_s(t)$ (red) at time $t=1.15$ as a solution of the sharp-interface model. The interfaces are visualized using thick white curves and the deformation in the solid phase is visualized using a thin white mesh.}
\label{fig:3d}
\end{figure}
}

Before we explain the diffuse-interface model and the sharp-interface model, we would like to make a comment on simplifying modeling assumptions and our approach.  
\renewcommand\labelitemi{\tiny$\bullet$}
\begin{itemize}[leftmargin=.25in]  
\item We neglect inertia to avoid effects caused by elastic or capillary waves.
\item We formulate the viscoelastic models entirely in Lagrangian coordinates.
\item  We assume no-slip conditions between phases and all phases have the same  viscosity $\mu$.
\item  We consider a multiphase system with solid, liquid and gas phases, where the liquid and gas phases are characterized by the absence of elastic properties. 
\item  We consider slightly compressible phases via $\tfrac{\kappa}{2}(\det\bF-1)^2$ in the elastic energy with $\kappa\gg 1$. Locking is avoided  by  using sufficiently high polynomial degree in the FE discretization \cite{auricchio2013approximation}.
\end{itemize}
While we made these choices to make it easier to follow the presentation or avoid certain technical difficulties, none of these choices 
should restrict the validity of our findings concerning the sharp-interface limit in the presence of moving contact lines.

\subsection{Diffuse-interface model}
First we present a ternary phase-field model, that is capable to describe the evolution of a liquid phase and a soft elastic phase surrounded by a gas phase. 
Therefore, within a given fixed box $\Omega\subset\mathbb{R}^d$ for $d=2,3$ let $\psi:[0,T]\times\Omega\to\mathbb{R}^2$ be a ternary phase-field and $\bchi_{\varepsilon}:[0,T]\times\Omega\to\mathbb{R}^d$ a deformation field. We further decompose the phase-field $\psi=(\psi_s,\psi_\ell)$ into solid and liquid phase, so that the remaining air phase $\psi_a$ is determined implicitly by $\psi_a=-1-\psi_s-\psi_\ell$ so that $\psi_a=+1$ if $\psi_s=\psi_\ell=-1$.
Using this convention, $\psi_i=1$ indicates the presence of phase $i$ and $\psi_i=-1$ the absence of phase $i$ for $i\in\{s,\ell,a\}$.  We denote the deformation gradient $\bF=\nabla\bchi_{\varepsilon}$ and the Jacobian determinant $J=\det(\bF)$. For $\qd=(\bchi_{\varepsilon},\psi)\in\mathcal{Q}_\varepsilon$, the cornerstone of the phase-field model is the free energy  $\mathscr{F}_{\varepsilon}:\mathcal{Q}_\varepsilon\to\mathbb{R}$  
\begin{align}\label{equ:diffuse_lagrange_energy}
    \mathscr{F}_{\varepsilon}(\qd) & = \int_{\Omega}W^{\varepsilon}[\qd]\dx, \qquad W^{\varepsilon}[\qd]=\Wel(\bF,\psi) + \Wpfe(\psi,\bF^{-T}\nabla\psi),
\end{align}
for which here we choose 
\begin{subequations}
\begin{align} 
    \label{equ:neo_hooke}   
    \Wel(\bF,\psi) 
    & = \frac{G(\psi)}{2}{\rm tr}(\bF^T\bF-I) + \tfrac{\kappa}{2}(J-1)^2\\
    \Wpfe(\psi,\bF^{-T}\nabla\psi) 
    &= \sum_{i\in\{s,\ell,a\}}\gamma_i\left[\tfrac{\varepsilon}{2}|\bF^{-T}\nabla\psi_i|^2+ \tfrac{1}{4\varepsilon}(1-\psi_i^2)^2\right]J
\end{align}
\end{subequations}
with surface parameters $\gamma_i$ that correspond to the usual surface tensions via
\begin{align}
\label{eqn:surftension}
    \gamma_{sa} = \frac{3}{2\sqrt{2}}(\gamma_s + \gamma_a),\quad
\gamma_{s\ell} = \frac{3}{2\sqrt{2}}(\gamma_s + \gamma_\ell),\quad
\gamma_{{\ell}a} = \frac{3}{2\sqrt{2}}(\gamma_\ell + \gamma_a)\,. 
\end{align}
The function $G(\psi)$ in \eqref{equ:neo_hooke} interpolates the elastic properties of the different phases by
\begin{align}
G(\psi)=\tfrac12 G_\ell (1+\psi_\ell)+\tfrac12 G_s (1+\psi_s)+\tfrac12 G_a (1+\psi_a)\,,
\end{align}
with $G_i\in\R$, but we are going to set $G_\ell=G_a=0$ and rescale $G_s=1$. The material is nearly incompressible, which we achieve by setting $\kappa\gg 1$. On the boundary of the domain we use either homogeneous Dirichlet boundary conditions $\bchi_{\varepsilon}(x)-x=0$ for the entire displacement or for the normal component. In the latter case, homogeneous natural boundary conditions are implied for the tangential component.

Motivated by the \emph{extended gradient structures} concept introduced in \cite{schmeller2022gradient}, we use a dynamical model for the evolution $\qd:[0,T]\to\mathcal{Q}_\varepsilon$ with chemical potentials $\eta:[0,T]\to\mathcal{U}$ that satisfy
\begin{align}
\label{eqn:extgrad_di}
    \begin{split} a_{\varepsilon}(\eta,\phi)+b(\phi,\partial_t{\qd}) &=0\,,\\
    b(\eta,v) 
    s(\partial_t \bchi_{\varepsilon}, \boldsymbol{v}_{\bchi})
    &= \langle {\rm D}\calF_{\varepsilon}(\qd),v\rangle_{\calVd}\,,
    \end{split}
\end{align}
for all $\phi\in\mathcal{U}=H^1(\Omega;\R^2)$ and all $v=(\boldsymbol{v}_{\bchi},v_\psi)\in\calVd=H^1(\Omega;\R^d)\times H^1(\Omega;\R^2)$ and given initial values $\qd(t=0)\in\mathcal{Q}_\varepsilon$. For this model we use the bilinear forms
\begin{subequations}
\begin{align} 
&a_\varepsilon(\eta,\phi)=\int_\Omega m(\varepsilon)\bF^{-T}\nabla\eta\cdot\bF^{-T}\nabla\phi\dx,\\
&b(\phi,v)=\int_\Omega \phi v_\psi\dx,\\
&s(\boldsymbol{w},\boldsymbol{v})=\int_\Omega \mu(\psi)\,((\nabla\boldsymbol{w})\bF^{-1}):((\nabla\boldsymbol{v})\bF^{-1})\det\bF\dx,
\end{align}
\end{subequations}
with phase-dependent viscosity function
\begin{align}
\mu(\psi)=\tfrac12 \mu_\ell (1+\psi_\ell)+\tfrac12 \mu_s (1+\psi_s)+\tfrac12 \mu_a (1+\psi_a)\,,
\end{align}
with $\mu_i\in\R$. For simplicity here we choose $\mu_i=\mu$, such that the viscosity is constant $\mu=\mu(\psi)$.
By testing \eqref{eqn:extgrad_di} with $v=\partial_t \qd$ and $\phi=\eta$ we directly obtain thermodynamic consistency
\begin{align}
\frac{\rm d}{{\rm d}t}\calF_{\varepsilon}
=\langle {\rm D}\calF_{\varepsilon}(\qd),\partial_t \qd\rangle_{\mathcal{V}} = - (s(\partial_t\bchi_{\varepsilon},\partial_t\bchi_{\varepsilon})+a_\varepsilon(\eta,\eta))\le 0.
\end{align}
Note that for clarity of the presentation, the dependence of the bilinear forms on the state variable $\qd$ is not shown explicitly. 
\begin{remark}
We make all the computations for $d=3$ by assuming axisymmetry. Therefore, the weak formulation in cylindrical coordinates is obtained by replacing the domain with a cylinder of radius $R$ and height $H$ such that $\Omega_r=[0,R]\times[0,H]$. Using $\boldsymbol{w}=(w_r,w_z)$ we make the following replacements
\begin{align*}
\dx = 2\pi\mathrm{d}r\,\mathrm{d}z,\quad
\nabla\boldsymbol{w} = \begin{pmatrix} \partial_r w_r & 0 & \partial_z w_r \\
0 & r^{-1}w_r & 0 \\
\partial_r w_z & 0 & \partial_z w_z
\end{pmatrix},\quad
\nabla \eta = \begin{pmatrix}
\partial_r \eta \\ 0 \\ \partial_z \eta \end{pmatrix},
\end{align*}
for the integration measure $\dx$, for gradients of vector fields $\nabla\boldsymbol{w}$ and for gradients of scalar fields $\nabla\eta$. These replacements are performed consistently in the bilinear forms and the free energy. 
\end{remark}

\subsection{Sharp-interface model}
Now we present the corresponding sharp-interface model of the multiphase model that is our candidate for the sharp-interface limits as $\varepsilon\to 0$. Using the same fixed box $\Omega\subset\mathbb{R}^d$ 
instead we consider that each phase $i\in\{s,\ell,a\}$ occupies a sufficiently regular subdomain $\Omega_i\subset\Omega$ that do not overlap $\Omega_i\cap\Omega_j=\emptyset$ and fill the entire domain, i.e., $\Omega =\overline{\Omega_s\cup\Omega_\ell\cup\Omega_a}$.
For $q_0\equiv \bchi_0 \in\mathcal{Q}_0$, the cornerstone of the sharp-interface model is also a free energy  $\mathscr{F}_{0}:\mathcal{Q}_0\to\mathbb{R}$  
\begin{align}\label{equ:sharp_lagrange_energy}
\calF_0(\qs)=\sum_{i\in\{s,\ell,a\}}\int_{\Omega_i}   \Wel^i(\bF) \dx + \sum_{ij\in\{s\ell,sa,\ell a\}}\int_{\Gamma_{ij}} \gamma_{ij}|\mathrm{cof}(\bF)\cdot\boldsymbol{\nu}|\,{\rm d}s\,,
\end{align}
where $\bF=\nabla\bchi_0$ as before and 
for which here we choose  
\begin{align} 
    \Wel^i(\bF) 
    & = \frac{G_i}{2}{\rm tr}(\bF^T\bF-I) + \tfrac{\kappa}{2}(J-1)^2,
\end{align}
and for the surface tensions $\gamma_{ij}$ from \eqref{eqn:surftension} associated to the (sharp) interfaces $\Gamma_{ij}=\partial\Omega_i\cap\partial\Omega_j$ for $ij\in\{s\ell,sa,\ell a\}$. Also motivated by a gradient structure, we use a dynamical model for the evolution $\qs:[0,T]\to\mathcal{Q}_0$ that satisfies
\begin{align}
\label{eqn:extgrad_si}
     -
    s(\partial_t \bchi_0, \boldsymbol{v}_{\bchi})
    &= \langle {\rm D}\calF_{0}(\qs),\boldsymbol{v}_{\bchi}\rangle_{\mathcal{V}}\,,
\end{align}
for all $\boldsymbol{v}_{\bchi}\in\mathcal{V}_0 =H^1(\Omega;\R^d)$ and given initial values $q_0(t=0)\in\mathcal{Q}_0$. For this model we use the same bilinear form as before, i.e.
\begin{align}
s(\boldsymbol{w},\boldsymbol{v})=\sum_{i\in\{s,\ell,a\}}\int_{\Omega_i} \mu_i\,((\nabla\boldsymbol{w})\bF^{-1}):((\nabla\boldsymbol{v})\bF^{-1})\det\bF\dx\,,
\end{align}
with phase-dependent viscosities $\mu_i\in\R$, which here we set equal to $\mu_i=\mu$.
By testing \Cref{eqn:extgrad_si} with $\partial_t \bchi_0$ we directly obtain thermodynamic consistency
\begin{align}
\frac{\rm d}{{\rm d}t}\calF_{0}
=\langle {\rm D}\calF_{0}(\qs),\partial_t \bchi_0\rangle_{\mathcal{V}_0} = - s(\partial_t\bchi_0,\partial_t\bchi_0)\le 0.
\end{align}

\begin{remark}
Assuming axisymmetry, we additionally need to replace
\begin{align*}
&\mathrm{d}s_x = 2\pi r\mathrm{d}s_r,
\quad \boldsymbol{\nu}=\begin{pmatrix}\nu_r \\ 0 \\ \nu_z \end{pmatrix},
\end{align*}
the surface measure $\mathrm{d}s_x$ and the normal vector $\boldsymbol{\nu}$ in the energy $\calF_0$.
\end{remark}

\begin{remark}[Eulerian formulation]
The Eulerian free energy for the phase-field model  \eqref{equ:diffuse_lagrange_energy} with $\bar{q}_\varepsilon=(\bar{\bF},\bar{\psi})$ is 
\begin{align}
    \calF_\varepsilon(\bar{q}_{\varepsilon})=\int_{\Omega}  \bar{J}^{-1}\left(\Wel(\bar\bF,\bar\psi) + \Wpfe(\bar\psi,\nabla\bar\psi)\right) \,{\rm d}\bar{x} \,,
\end{align}
where $\bar{\psi}\circ\bchi_\varepsilon=\psi$ and $\bar\bF\circ\bchi_\varepsilon=\bF$ and $\bar{J}=\det\bar\bF$. 
The free energy of the sharp-interface model \eqref{equ:sharp_lagrange_energy} in the Eulerian frame is
\begin{align}
\calF_0(\bar{q}_0)=\sum_{i\in\{s,{\ell},a\}}\int_{\bar\Omega_i}   \bar{J}^{-1}\Wel^i(\bar\bF) \,{\rm d}\bar{x} + \sum_{ij\in\{s{\ell},{\ell}a,sa\}}\int_{\bar\Gamma_{ij}} \gamma_{ij}\,{\rm d}\bar{s}\,, 
\end{align}
where the Eulerian state $\bar{q}_0=(\bar\bF,\bar\Omega_i)$ contains  $\bar\Omega_i(t)=\bchi_0(t,\Omega_i(0))$ the domain currently occupied by phase $i\in\{s,\ell,a\}$ at time $t$. The PDEs in Eulerian variables also contain a velocity $\bar\bu:[0,T]\times\Omega\to\mathbb{R}^d$, which, however, the free energy does not depend on for viscous flows without kinetic energy.
\end{remark}

\EEE
\subsection{Scaling limits for the Cahn-Hilliard mobility}\label{sec:scaling_limits}

Varying choices/scalings of the Cahn-Hilliard mobility in $\varepsilon$, may approximate different sharp-interface models.
Generally speaking, a bound $\mob(\varepsilon)\lesssim\mathcal{O}(\varepsilon)$ on the mobility is needed to prevent artificial diffusion of the phase fields \cite{huang2009mobility}. 
A further bound from the other side is essential to ensure that the mobility does not become too small, which would lead to a different/unintended sharp-interface limit \cite{schaubeck2014sharp}.

{A fundamental work on the passage from the phase-field model to the sharp-interface limit is \cite{abels2012thermodynamically}. The sharp interface is analyzed using matched asymptotic techniques, and degenerate mobilities are included in their study. 
Precisely, the following four cases are discussed:  
\begin{align*}
    m({\varepsilon},\psi)
    = 
    \begin{cases}
    m_0 & \text{Case 1}\\
    \varepsilon m_0 & \text{Case 2}\\
    \frac{m_1}{\varepsilon}(1 - \psi^2)_+ & \text{Case 3}\\
    m_1(1 - \psi^2)_+ & \text{Case 4} 
    \end{cases}\,,
\end{align*}
where for Case 2 and 4 a sharp-interface limit corresponding to ours is reached. In the other cases, additional terms appear in the kinematic conditions of the sharp-interface model. 

Having introduced the diffuse and the associated sharp-interface model, we will now analyze the limit passage $\varepsilon\to 0$.
In this work we want to generalize some considerations concerning viable sharp-interface limits, for which the passage from the diffuse-interface model \eqref{equ:diffuse_lagrange_energy} to the sharp-interface model \eqref{equ:sharp_lagrange_energy} is valid including contact lines. The sharp-interface limit, $\varepsilon\to 0$, of the diffuse-interface model depends substantially on the scaling of the Cahn-Hilliard mobility $\mob=m_0\varepsilon^{\alpha}$, $0\leq\alpha\leq\infty$. 
Degenerate mobilities or other forms of Cahn-Hilliard mobility are also used and require appropriate scaling properties \cite{abels2012thermodynamically,lee2016sharp,dziwnik2017anisotropic}.
\EEE

We believe the setup shown in \Cref{fig:3d_drops} has several mathematical and numerical challenges that make it a suitable benchmark to study the sharp-interface limit of phase-field models: 
\begin{itemize}[leftmargin=.25in]
\item  For this setup, the stationary state is close to the initial data so that Lagrangian methods can compute a regular deformation $\bchi$ without needing ALE methods or remeshing.
\item  With the initially flat solid substrate, surface energies will quickly enforce a contact angle (kink in the solid surface) via the Neumann triangle construction. While the resulting nonsmoothness in $\bF$ is a challenge, this scenario is typical for soft wetting applications. This could lead to suboptimal convergence rates in space discretizations.

\item  If initial data do not satisfy the Neumann triangle construction, then solutions are also expected to show nonsmooth behavior in time, leading to suboptimal convergence in time.
\end{itemize}

\begin{figure}[hb!]
    \centering
a)\,\includegraphics[height=0.32\textwidth,trim=1.5cm 0cm 1.5cm 0cm,clip]{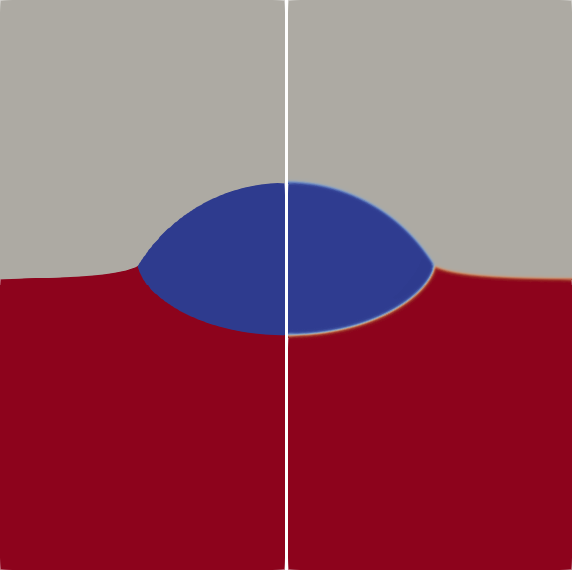}\,\,\,%
b)\,\includegraphics[height=0.32\textwidth,trim=1.5cm 0cm 1.5cm 0cm,clip]{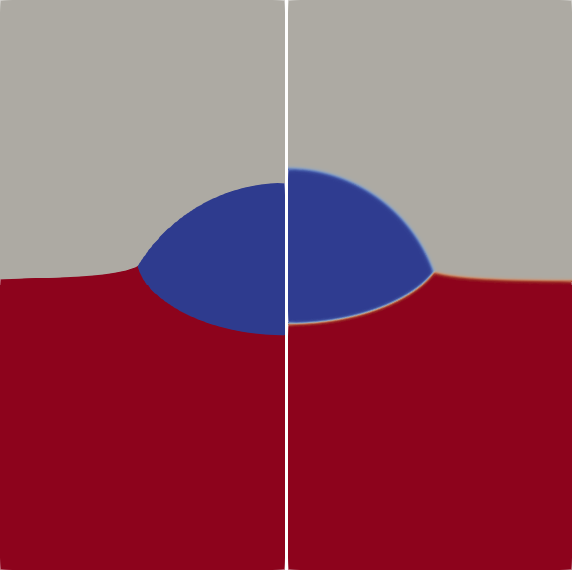}\,\,\,%
c)\includegraphics[height=0.32\textwidth,trim=10.5cm 12cm 17cm 9cm,clip]{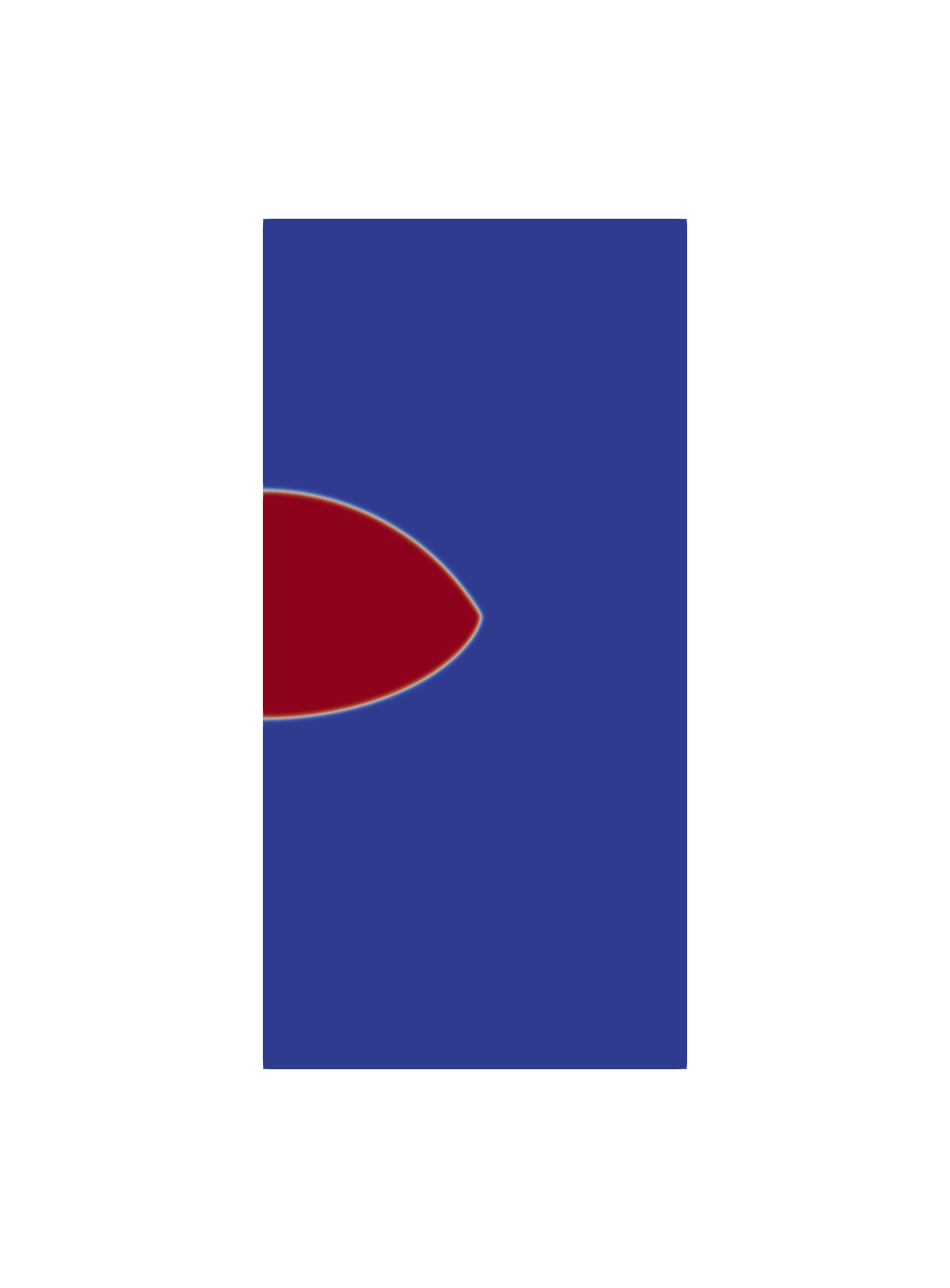}\!%
\includegraphics[height=0.32\textwidth,trim=10.5cm 12cm  6.7cm 9cm,clip]{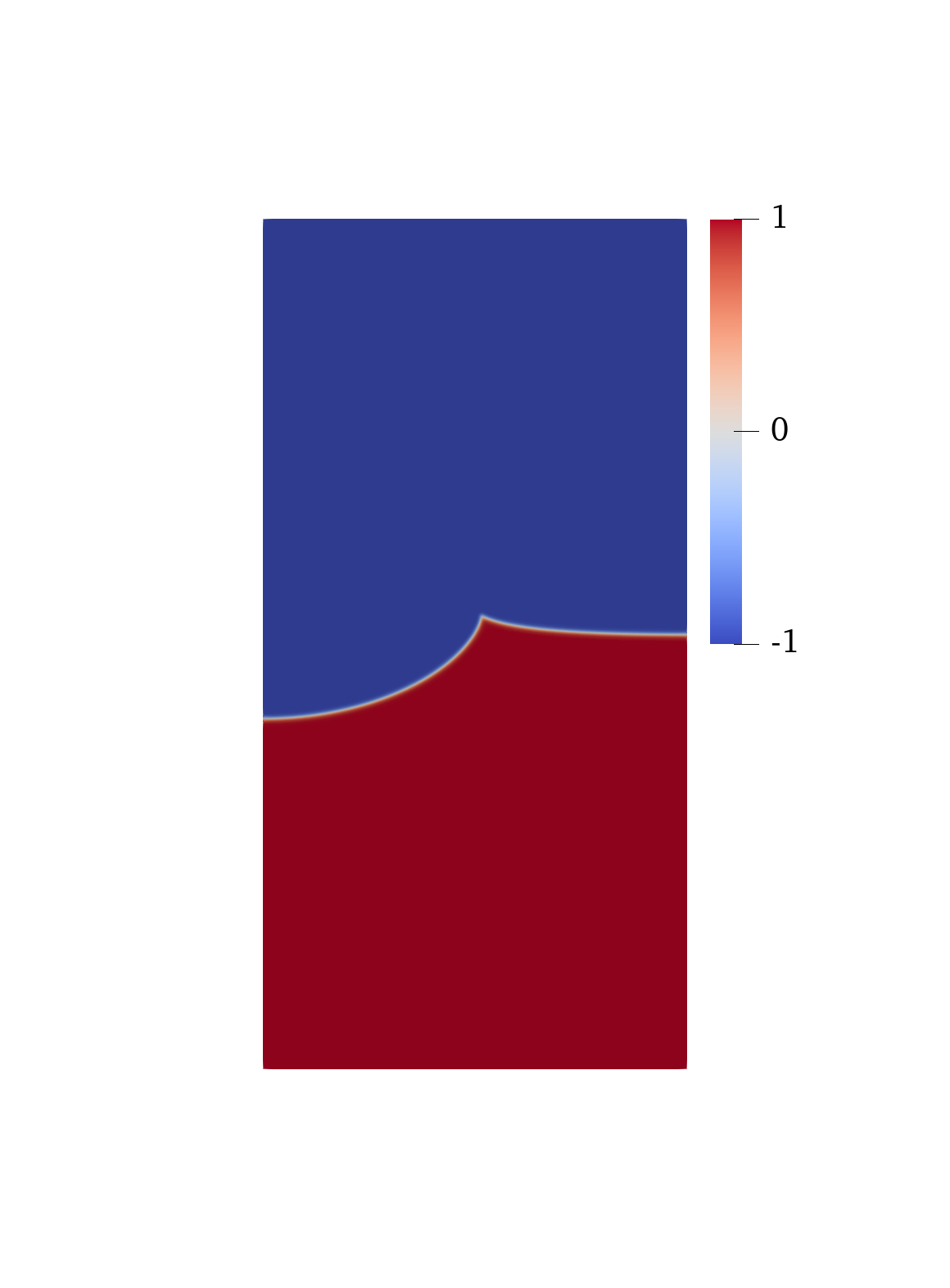}%
    \caption{In a,b) the left side shows the sharp interface at time $t=T=1.15$ and the right side is the diffuse interface model for mobility a) $m=\varepsilon^{5/2}$ and b) $m=0$ with $\varepsilon=2\bar{\varepsilon}$ in both cases. The deformed phase fields are plotted via $\bar{\psi}_{\rm id}=(2+\bar{\psi}_s-\bar{\psi}_\ell)/2$ taking values $\bar{\psi}_\mathrm{id}=0,1,2$ in the liquid  (blue), air (gray), solid phase (red), respectively.
    c) Deformed diffuse phase fields for a) showing (left) liquid phase $-1\le\bar{\psi}_\ell\le 1$ and (right) solid phase $-1\le\bar{\psi}_s\le 1$.}
    \label{fig:mob_to_small}
\end{figure}

In the following, in order to improve the visibility of results we will show solutions of the sharp-interface model and of the phase-field model only as cross sections in the $(r,z)$-plane.
\Cref{fig:mob_to_small} shows in the left half of each droplet the result of the sharp-interface at fixed time $t=1.15$, while the right half shows the diffuse-interface model with $\varepsilon = 2\bar{\varepsilon}$ with $\bar{\varepsilon}=1.8775\cdot 10^{-3}$ at the same time with two different mobilities. For $\mob=\varepsilon^{\alpha}$ with $\alpha = 5/2$ used in the left/first figure, we see a comparison. 
In the second/right droplet we choose $\mob = 0$ and a clear deviation between the diffuse and the sharp-interface model can be seen. 
{The reason for this deviation is that the phase-field $\psi$ upon deformation to $\bar\psi\circ\bchi_\varepsilon=\psi$ can not relax to its optimal $\tanh$-profile and therefore the approximation of the Eulerian surface energy is not guaranteed \cite{schaubeck2014sharp}.}
Thus, this indicates that the Cahn-Hilliard mobility is chosen too small, i.e. $\alpha<\underline{\alpha}$.
Furthermore, the upper bound to the Cahn-Hilliard mobility must also be respected. 
The next step is to numerically analyze the convergence of the models in appropriate norms.

\section{Discretization in space and time}
\label{sec:si_vs_di}
In this section, we numerically estimate the distance of the the phase-field and sharp-interface models to make a statement about the convergence to the sharp-interface limit. Therefore, the two main assumptions that need to be verified are:\\[-0.5em]


\begin{enumerate}
\item[\bf A1]  The deformation converges $\|\bchi_0-\bchi_\varepsilon\|\to 0$ as $\varepsilon\to 0$ for a suitable norm $\|\cdot\|$.

\item[\bf A2]  The (Lagrangian) phase fields $\psi_i$ remain close to $\Omega_i$ in the sense
\begin{align*}
\mathrm{distance}(\Omega_{i,\varepsilon}(t),\Omega_i)\to 0,
\end{align*}
as $\varepsilon\to 0$ for time $t$ and a suitable distance between $\Omega_{i,\varepsilon}(t)=\{{x}\in\Omega:\psi_i(t,{x})>0\}$ and the (time-independent) sharp domains $\Omega_i$ for $i\in\{s,\ell,a\}$. 
\end{enumerate}
\vspace{0.2em}
\noindent
In order to make meaningful statements about this convergence, we need to estimate the numerical errors in the respective distances and norms. Let us denote by $\bchi_{\varepsilon,\Delta}$ a numerical solution of the diffuse-interface model for fixed $\varepsilon$ and by 
$\bchi_{0,\Delta}$ a numerical solution of the sharp-interface model. As above by $\bchi_\varepsilon$ and $\bchi_0$ we denote the exact solutions. 
Then, using the triangle inequality, we have the estimate for the sharp-interface limit
\begin{align}
\underbrace{\|\bchi_0-\bchi_\varepsilon\|}_{\text{sharp-interface limit}} \le \underbrace{\|\bchi_0 - \bchi_{0,\Delta}\|}_{\text{error sharp interface $e_{0,\Delta}$}} + \quad\underbrace{\|\bchi_{0,\Delta} - \bchi_{\varepsilon,\Delta}\|}_{\text{estimate sharp-interface limit}}\quad+ \underbrace{\|\bchi_{\varepsilon,\Delta} - \bchi_\varepsilon\|}_{\text{error diffuse interface $e_{\varepsilon,\Delta}$}}
\end{align}
in an appropriate function space norm $\|\cdot\|$. The first goal here is to estimate various contributions to the approximation errors $e_{0,\Delta}=\|\bchi_0 - \bchi_{0,\Delta}\|$ and $e_{\varepsilon,\Delta}=\|\bchi_{\varepsilon,\Delta} - \bchi_\varepsilon\|$ for simulation parameters associated with the discretization such as $\Delta=\{\tau,h,k_{\bchi},k_{\psi}\}$ for time-step size $\tau$, mesh size $h$, polynomial degree for deformation $k_{\bchi}$, and for the polynomial degree for the phase fields $k_\psi$. These we verify by time-step bisection $\tau\to\tau/2$, uniform refinement $h\to h/2$, and by varying the polynomial degree. Then we compute $\|e_{0,\Delta-\Delta'}\|=\|\bchi_{0,\Delta}-\bchi_{0,\Delta'}\|$ for different $\tau$ and $\Delta=\{\tau,h,k_{\bchi},k_{\psi}\}$ and $\Delta'=\{\tau/2,h,k_{\bchi},k_{\psi}\}$ to approximate the convergence of the time-discretization and for the other cases respectively. In such a case we write for simplicity $\|e_{0,\Delta-\Delta'}\|\equiv\|e_{0,\tau-\tau/2}\|$ and assume all other discretization parameters are known an fixed. For the space discretizations the computation of the error might also involve a projection $P_h:\mathcal{V}_{0,h}\to\mathcal{V}_{\varepsilon,h'}$ or 
$P_h:\mathcal{V}_{0,h}\to\mathcal{V}_{0,h'}$ or $P_h:\mathcal{V}_{\varepsilon,h}\to\mathcal{V}_{\varepsilon,h'}$, but this resulting projection error was mostly negligible compared to the error of the other approximations.

The error between the deformation fields is measured in $L^2$ and $L^{\infty}$ Bochner norms which is a concept that extends the classical Sobolev norms \cite{adams2003sobolev} to time-dependent problems via
\begin{align}
   \|\bchi\|_{L^2(Q_T;\R^d)}
   = \left(\int_0^T\|\bchi(t)\|^2_{L^2(\Omega;\mathbb{R}^d)}{\rm d}t\right)^{\tfrac{1}{2}}\,,\quad
   \|\bchi\|_{L^{\infty}(Q_T;\R^d)}
   = \esssup_{t\in[0,T]}\|\bchi(t)\|_{L^{\infty}(\Omega;\mathbb{R}^d)}\,.
\end{align}
While the $L^2$ Bocher norm measures the error over the whole area and neglects errors in small regions, the $L^\infty$ norm is also sensitive to errors in small regions, e.g. on interfaces and at contact lines.
First we introduce separately the space and time discretizations of phase-field and sharp-interface model.

\begin{remark}\label{rem:comment_to_A2} Using the deformations, we can map A2 to an Eulerian version and require
\begin{align}
\label{reform}
\mathrm{distance}(\bar\Omega_{i,\varepsilon}(t),\bar\Omega_i(t))\to 0,
\end{align}
as $\varepsilon\to 0$ for any time $t$ with $\bar\Omega_{i,\varepsilon}(t)=\bchi_\varepsilon(t,\Omega_{i,\varepsilon}(t))$ and $\bar\Omega_i(t)=\bchi_0(\Omega_i)$ and $i\in\{s,\ell,a\}$. Assuming convergence A1 in a sufficiently strong norm, then  A2  and \eqref{reform} should become equivalent. We ensure A2 by making $m(\varepsilon)$ small enough as $\varepsilon\to 0$.
\end{remark}

\subsection{Space and time discretization}

\paragraph*{Space-discretization} Based on the weak formulation of the diffuse-interface model \eqref{eqn:extgrad_di} and the sharp-interface model \eqref{eqn:extgrad_si} we employ the finite element method to derive a discretization in space. The main idea of this \emph{structure preserving} discretization is to adopt a weak formulation for the phase fields with $\mathcal{Q}_{h,\varepsilon}\subset\mathcal{Q}_\varepsilon$, $\mathcal{U}_h\subset\mathcal{U}$, $\mathcal{V}_{\varepsilon,h}\subset\mathcal{V}_\varepsilon$ and for the sharp-interface limit with 
$\mathcal{Q}_{h,0}\subset\mathcal{Q}_0$, $\mathcal{V}_{0,h}\subset\mathcal{V}_0$ via a finite element method. 

For the sharp-interface model we use computational meshes, where the elements and edges are aligned with the phases $\Omega_i$ and the interfaces $\Gamma_{ij}$ as shown in the first (left) panel of \Cref{fig:initial_meshes}. For the convergence study we usually use finer meshes for the sharp-interface model, where we perform 1-3 uniform refinements of the first mesh (mesh a) shown in \Cref{fig:initial_meshes} and project the vertices of $\Gamma_{{\ell}a}$ (interface between blue and gray domain) back onto the set $\sqrt{r^2+(z-1)^2}=\nicefrac12$. This coarsest mesh has 453 vertices resulting in the dimension $\dim\mathcal{V}_{h,0}=3{\,}496$, while after three uniform refinements the mesh has 27{\,}221 vertices resulting in the dimension $\dim\mathcal{V}_{h,0}=216{\,}786$ using $P_2$ finite elements for the deformation.

\begin{figure}[hb!]\label{fig:different_inital_meshes}
    \centering
    \includegraphics[width=0.66\textwidth]{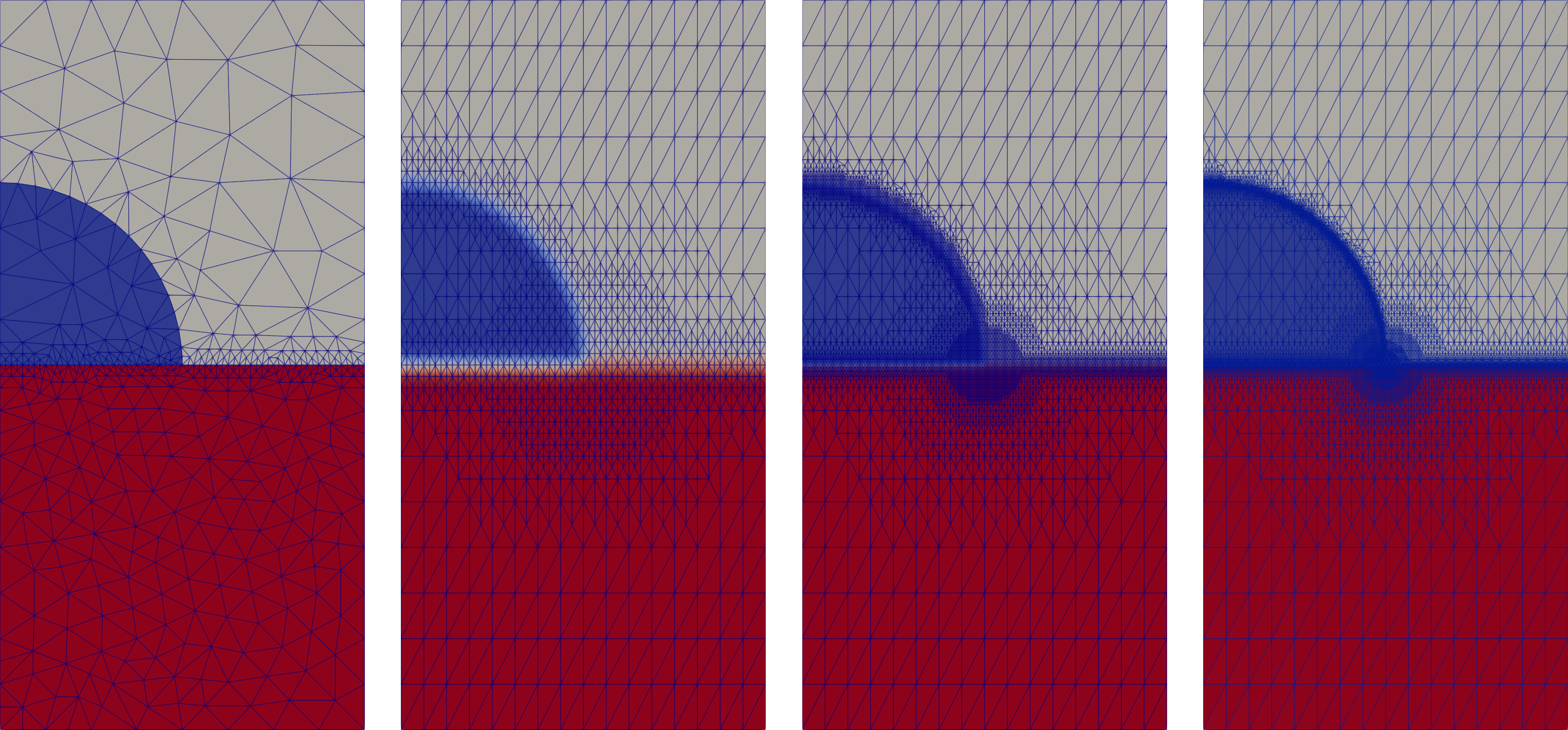}

    \includegraphics[width=0.66\textwidth]{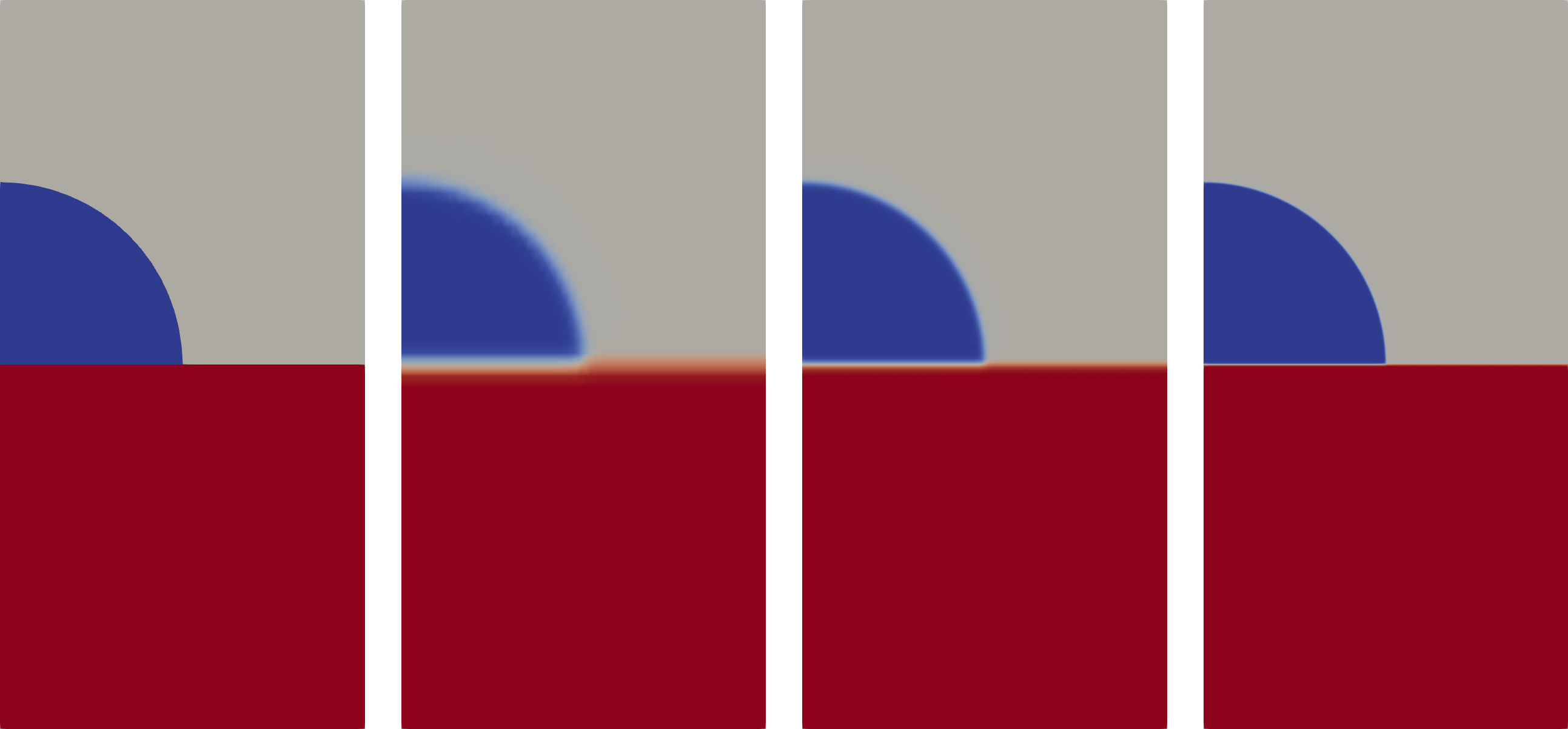}

     a) \hspace{2cm} b) \hspace{2cm} c) \hspace{2cm} d)
    
    \caption{(Top) Different initial meshes from (left) to (right) for a) sharp interface, b) phase-field with $\varepsilon=16\bar{\varepsilon}$, c) phase-field with $\varepsilon=4\bar{\varepsilon}$ and d) phase-field with $\varepsilon=\bar{\varepsilon}$ for $\bar{\varepsilon}=0.001875$ and (bottom) corresponding phase indicators. For the phase field we show $\psi_\mathrm{id}$.}
    \label{fig:initial_meshes}
\end{figure}

For the diffuse-interface model the base mesh strongly depends on the values of the interfacial thickness $\varepsilon$, where we consider for $\bar{\varepsilon}=0.001875$ the values $\varepsilon=n\bar{\varepsilon}$ for $n=1,2,4,8,16,32$. The corresponding meshes shown in \Cref{fig:initial_meshes} have $1{\,}310$ (mesh b), $5{\,}205$ (mesh c), $19{\,}502\}$ (mesh d) vertices corresponding to $\dim\mathcal{V}_{h,\varepsilon}=10{\,}322$, $\dim\mathcal{V}_{h,\varepsilon}=41{\,}458$, $\dim\mathcal{V}_{h,\varepsilon}=155{\,}810$ using $P_2$ finite elements for the deformation, respectively. On top of that, the phase-field model contains scalar unknowns for the phase fields $\psi$ and their chemical potentials $\eta$. Note that the mesh is refined near the interfaces to resolve the transition layer of width $\varepsilon$ and in a larger region near the contact line to account for potential artificial diffusion.

Let $P_k(T,\R^r)$ be the set of $\R^r$-valued polynomials of degree $k$ restricted to triangles $T\in\mathcal{T}_h$ and $\smash{\Omega=\bigcup_{T\in\mathcal{T}_h} T}$ an admissible triangulation of the ternary system. We use spaces  $\smash{\boldsymbol{V}^{k}}$ for vectorial unknowns and $V^{k}$ for scalar unknowns, where we define the $H^1$-conforming finite element spaces
\begin{align*}
&\boldsymbol{V}_h^{k}=\{\boldsymbol{v}\in C^1(\Omega,\R^d): \boldsymbol{v}|_T\in P_{k}(T,\R^d), T\in\mathcal{T}_h\},\\
&V_h^{k}=\{v\in C^1(\Omega,\mathbb{R}): v|_T\in P_{k}(T,\R), T\in\mathcal{T}_h\}\,.
\end{align*}
For the sharp-interface model we use $\mathcal{V}_{h,0}=\boldsymbol{V}_h^{k_{\bchi}}$ with $k_{\bchi}=1,2,3$
and enforce essential boundary conditions when solving the corresponding nonlinear problem and for the diffuse-interface model we use 
\begin{align*}
\mathcal{V}_{h,\varepsilon}=\boldsymbol{V}_h^{k_{\bchi}}\times V_h^{k_\psi}\times V_h^{k_\psi}, \qquad \mathcal{U}_h=V_h^{k_\psi}\times V_h^{k_\psi},
\end{align*}
with $k_{\bchi}=1,2,3$ and $k_\psi=2,3$. As for the sharp-interface model, essential boundary conditions for the deformation are enforced when solving the corresponding nonlinear problem. Nonlinear elasticity problems with compressibility $\kappa\gg 1$ are likely to show locking phenomena, which we deal with by choosing $k_{\bchi}$ large enough.\\[-0.7em]

\paragraph*{Time-discretization}
The time discretization is constructed as in \cite{schmeller2022gradient} based on an incremental-minimization-scheme for $q_\varepsilon^k=q_\varepsilon(k\tau)$ and $q_0^k=q_0(k\tau)$ with time-step size $\tau$. Based on the previous formal definition of the weak formulations for the diffuse-interface model in \eqref{eqn:extgrad_di} and sharp-interface model in \eqref{eqn:extgrad_si} from  \Cref{sec:models} this leads to the following incremental nonlinear problems.\\[-0.7em]

\noindent
\textbf{Incremental scheme for sharp-interface model:} For given $q_0^{k-1}\equiv\bchi_0^{k-1}\in\mathcal{V}_{h,0}$ find $q^{k}_0\equiv\bchi_0^{k}\in\mathcal{V}_{h,0}$ such that
\begin{align}
\label{eqn:num_extgrad_si}
     -
    \tfrac{1}{\tau} s(\bchi^k_0-\bchi^{k-1}_0, \boldsymbol{v}_{\bchi})
    &= \langle {\rm D}\calF_{0}(\bchi^k_0),\boldsymbol{v}_{\bchi}\rangle_{\mathcal{V}_{h,0}}\,,
\end{align}
for all $\boldsymbol{v}_{\bchi}\in\mathcal{V}_{h,0}$. While the left-hand side of the problem appears is linear if in $s$ we use the previous state $\qs^{k-1}$, it is usually nonlinear through the dependence of $\mathrm{D}\calF_0(\bchi_0^k)$ on $\bchi_0^k$.\\[-0.7em]

\noindent
\textbf{Incremental scheme for diffuse-interface model:} 
For given $q_{\varepsilon}^{k-1}=(\bchi_{\varepsilon}^{k-1},\psi^{k-1})\in\mathcal{V}_{h,\varepsilon}$ and $\eta^{k-1}\in\mathcal{U}_h$ find $q_{\varepsilon}^k=(\bchi_{\varepsilon}^{k},\psi^k)\in\mathcal{V}_{h,\varepsilon}$ and $\eta^k\in\mathcal{U}_h$ such that
\begin{align}
\label{eqn:num_extgrad_di}
    \begin{split} a_{\varepsilon}(\eta^k,\phi)+\tfrac{1}{\tau}b(\phi,{\qd^k-\qd^{k-1}}) &=0\,,\\
    b(\eta^k,v) - 
    \tfrac{1}{\tau} s(\bchi^k_{\varepsilon}-\bchi^{k-1}_{\varepsilon}, \boldsymbol{v}_{\bchi})
    &= \langle {\rm D}\calF_{\varepsilon}(\qd^k),v\rangle_{\mathcal{V}_{h,\varepsilon}}\,,
    \end{split}
\end{align}
for all $v=(\boldsymbol{v}_{\bchi},v_\psi)\in\mathcal{V}_{h,\varepsilon}$ and all $\phi\in\mathcal{U}_h$. \\[-0.7em]

To initialize the chemical potentials $\eta(t=0)$ in the diffuse-interface model and to overcome any initial transient (in both models) we perform the first iteration with a small time-step size $0<\tau_0\ll \tau$. The next 20 iterations are performed with a time-step size $\tau/10$ and then we output the solutions $q_{0,\Delta}^k$ and $q_{\varepsilon,\Delta}^k$ every multiples of $0.05$ time units until the final time $T$ is reached.
Next we will discuss space-time discretization errors for numerical solutions 
\begin{align}
q_{0,\Delta} & \qquad\rightarrow\qquad \Delta = \{\tau,h,\boldsymbol{P}_{k_{\bchi}}\},\\
q_{\varepsilon,\Delta} & \qquad\rightarrow\qquad \Delta = \{\tau,h,\boldsymbol{P}_{k_{\bchi}},P_{k_\psi}\},
\end{align}
of the sharp-interface model $q_{0,\Delta}$ and the diffuse-interface model $q_{\varepsilon,\Delta}$ with the discretization parameters mentioned before, i.e. $h$ indicates the space discretization, $\tau$ is the time step size, $k_{\bchi}$ and $k_{\psi}$ denote the polynomial degree of the deformation and the phase fields, respectively.

\subsection{Estimation of numerical errors}
\label{sec:num_errors}

\paragraph*{Parameters}
The physical parameters for the sharp-interface and diffuse-interface model are listed in \Cref{tab:DI_parameters} and the corresponding surface tensions for the sharp-interface model can be derived from them using \Cref{eqn:surftension}. For integer $n\ge 0$, interfacial widths for the diffuse-interface model are  $\varepsilon=2^n\bar{\varepsilon}$ with  
$\bar{\varepsilon}=1.8775\cdot 10^{-3}$.
Our standard numerical parameters, i.e. the once causing the error collected in $\Delta$, are 
\begin{align}
    \tau = 0.005\,,\qquad k_{\bchi}  = 2\,,\qquad k_{\boldsymbol{\psi}} = 2\,,\qquad m_0=1\,,
\end{align}
and $\tau_0=0.5\cdot 10^{-7}$. The common mesh for the diffuse-interface model corresponds to $4\bar{\varepsilon}$ which is Case c) and for the phase-field model two uniform refinements of Case a) shown in  \Cref{fig:different_inital_meshes}. The time evolution of the energy\footnote{The energy is computed using the sharp-interface model with P$_2$ elements on the coarsest mesh a) from \Cref{fig:initial_meshes}.} is shown in \Cref{fig:energy}. We will consider solutions for $0<t<T=1.15$,  shown by the horizontal dotted line, after which the evolution can be considered stationary and no nontrivial contribution to the sharp-interface limit tested using the Bochner norms is expected. The energies of the two models agree well for $m=\varepsilon^{5/2}$ but are systematically lower for $m=\varepsilon^1$.\\[-0.7em]

\noindent
\textbf{Note:} We omitted the factor $2\pi$ in all surface and volume integrals, which modifies the $L^2$ error norms and energies correspondingly.

 \begin{figure}[hb!]
\centering
{\includegraphics[width=0.5\textwidth,trim=0cm 0.5cm 1cm 1.6cm,clip]{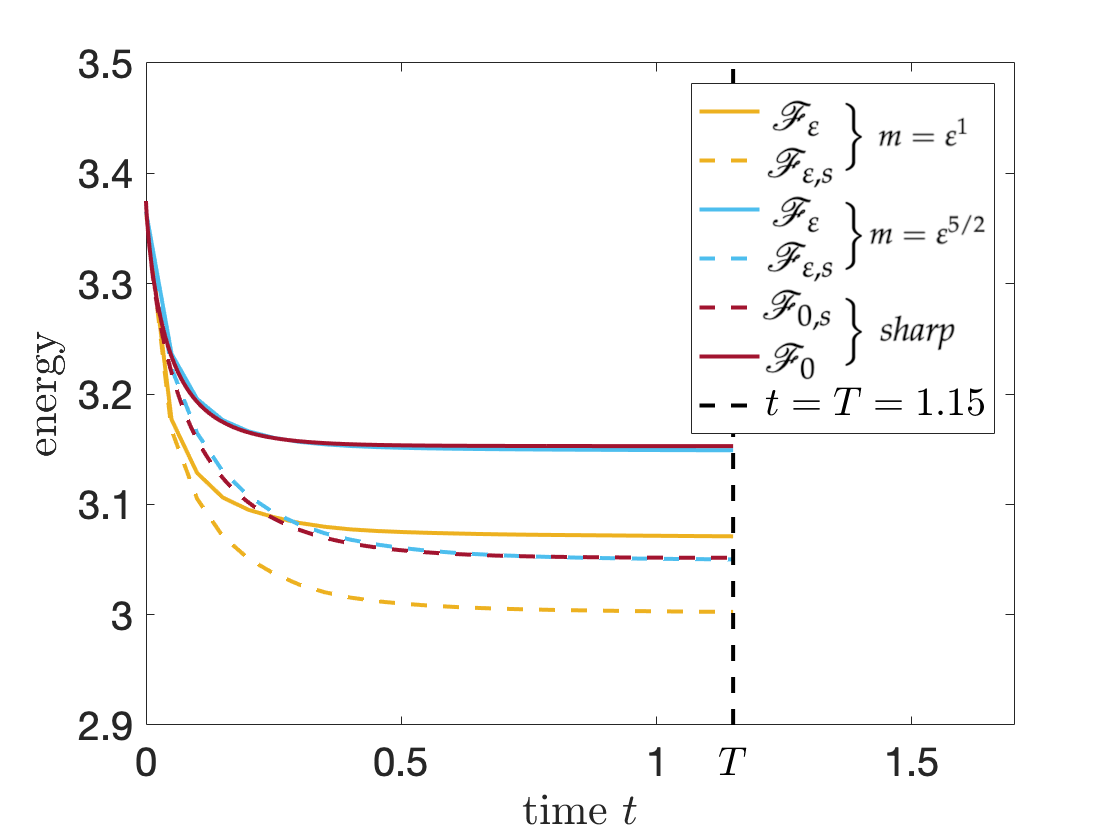}}
\caption{Free energies $\mathscr{F}_0$ and $\mathscr{F}_{\varepsilon}$ (full lines) and surface energies $\mathscr{F}_\textrm{0,s}=\sum_{ij}\int_{\bar{\Gamma}_{ij}}\gamma_{ij}\mathrm{d}\bar{s}$ and 
$\mathscr{F}_\textrm{$\varepsilon$,s}=\int_\Omega \Wpfe\dx$ (dashed lines) as a function of time $0\le t\le T$ for sharp (red) and diffuse interfaces  with mobilities $m=\varepsilon^1$ (yellow) and $m=\varepsilon^{5/2}$ (blue).}
\label{fig:energy}
\end{figure}

\begin{table*}[h]
\centering
\renewcommand{\arraystretch}{1.5}
  \begin{tabular}{cccccccccccc}
    \hline\hline
    parameter &$G_s$ & $G_{\ell,a}$ & $\gamma_s$ & $\gamma_\ell$ & $\gamma_a$ & $R$ &  $H$ & $r_\text{drop}$  &  $\kappa$ & $T$ \\
    \hline
    value & $10$ & $0$ & $0.5$  & $2.5$ & $3.5$ & $1.0$ & $2.0$ & $\nicefrac12$  & $10^{4}$ & $1.15$\\
    \hline\hline\\
  \end{tabular}
\caption{Parameters for the sharp-interface and diffuse-interface model: $G_i$ the shear modulus in the elastic energy, $\gamma_i$ the surface parameter for $i\in\{s,\ell,a\}$,  $R$ is the width of the rotational domain, $r_\text{drop}$ is the initial radius of the droplet, $h$ is the height of the substrate, $H$ is the height of the domain, and $T$ the final time.}
\label{tab:DI_parameters}
\end{table*}

\paragraph{Diffuse-interface model}
The numerical errors of the phase-field model measured using the two Bochner norms are shown in \Cref{tab:num_error_DI}. 
We observe clear convergence trends in space and time. As expected, errors in the $L^2$ norm are always smaller then those in the corresponding $L^\infty$ norm. The significant decrease in error between $e_{\varepsilon,h-\nicefrac{h}{2}}$ and $e_{\varepsilon,\nicefrac{h}{2}-\nicefrac{h}{4}}$ by a factor ten is  unclear but gives a comparable order as the time-discretization. For the time-discretization, we observe a reduction of the error between $e_{\varepsilon,\tau-2\tau}$ and $e_{\varepsilon,\tau-\nicefrac{\tau}{2}}$ by roughly a factor of $2$, which indicates the expected linear convergence rate in $\tau$.
By far the largest errors on the order of $10^{-3}$ for the $L^\infty$ norm emerge from the polynomial order, which indicate possible locking effects for $k_{\bchi}$ and a limiting resolution a the interface for $k_\psi$.

\renewcommand{\arraystretch}{1.6}
\begin{table}[H]
\centering
\begin{tabular}{c|c|c|c|c|c|c}
    norm    & $e_{\varepsilon,h-\nicefrac{h}{2}}$ & $e_{\varepsilon,\nicefrac{h}{2}-\nicefrac{h}{4}}$ & $e_{\varepsilon,\tau-2\tau}$ & $e_{\varepsilon,\tau-\nicefrac{\tau}{2}}$ & $e_{\varepsilon,\psi_{P2}-\psi_{P3}}$ & $e_{\varepsilon,\bchi_{P2}-\bchi_{P3}}$\\\hline
     $L^{2}$
     & 
     $2.8\cdot 10^{-4}$ & $2.7\cdot 10^{-5}$ & $9.5\cdot 10^{-5}$ & $5.3\cdot 10^{-5}$ & $1.2\cdot 10^{-4}$ & $2.9\cdot 10^{-4}$\\
     $L^{\infty}$
     & 
     $4.6\cdot 10^{-3}$ &  $3.2\cdot 10^{-4}$ & $9.5\cdot 10^{-4}$ & $5.5\cdot 10^{-4}$ & $2.7\cdot 10^{-3}$  & $5.0\cdot 10^{-3}$ 
\end{tabular}
\caption{Numerical errors for the phase-field model with respect to the $L^{2}(Q_T;\R^2)$ norm and with respect to the $L^{\infty}(Q_T;\R^2)$ norm for vector-valued and scalar functions upon change of time-step size $\{2\tau,\tau,\tau/2\}$, uniform refinement $\{h,\nicefrac{h}{2},\nicefrac{h}{4}\}$, and polynomial order $\{P_2,P_3\}$ of the FE-space of $\psi$ and $\bchi_\varepsilon$.}
\label{tab:num_error_DI}
\end{table}

\paragraph{Sharp-interface model}
For the sharp-interface model we show the numerical errors in \Cref{tab:num_error_SI}. For the uniform refinement we see a reduction in errors with a factor between $2$ and $4$ going from $e_{0,\nicefrac{h}{2}-\nicefrac{h}{4}}$ to $e_{0,\nicefrac{h}{4}-\nicefrac{h}{8}}$, which indicates a convergence order between 1 and 2 in space, as expected for $P_2$ elements and a less regular solution. The reduction of the errors $e_{0,\tau-\nicefrac{\tau}{2}}$ and $e_{0,\nicefrac\tau{2}-\nicefrac{\tau}{4}}$ is by a factor of two, showing linear convergence. As for the diffuse-interface model, by far the largest error come from different polynomial orders for the displacement and are of order $10^{-3}$. This is again indicative of possible locking effects.

Both diffuse-interface model and sharp-interface model clearly show numerical convergence, with the largest contribution to the error coming from polynomial degree and being of the order $10^{-3}$ in the $L^\infty$ norm.

\renewcommand{\arraystretch}{1.6}
\begin{table}[H]
\centering
\begin{tabular}{c|c|c|c|c|c|c|c}
    norm    & $e_{0,h-\nicefrac{h}{2}}$ & $e_{0,\nicefrac{h}{2}-\nicefrac{h}{4}}$ & 
    $e_{0,\nicefrac{h}{4}-\nicefrac{h}{8}}$ & $e_{0,\tau-\nicefrac{\tau}{2}}$ & $e_{0,\nicefrac{\tau}{2}-\nicefrac{\tau}{4}}$ & $e_{0,{P_1}-{P_2}}$ & $e_{0,{P_2}-{P_3}}$\\\hline
    $L^{2}$
    & $1.4\cdot 10^{-3}$ & $2.8\cdot 10^{-4}$ & $7.9\cdot 10^{-5}$ & $5.3\cdot 10^{-5}$ & $2.8\cdot 10^{-5}$ &
  $9.7\cdot 10^{-3}$ & $8.9\cdot 10^{-5}$\\
  $L^{\infty}$
  & $1.3\cdot 10^{-2}$ & $4.9\cdot 10^{-3}$ & $3.7\cdot 10^{-3}$ & $5.5\cdot 10^{-4}$ & $2.9\cdot 10^{-4}$ &
  $6.8\cdot 10^{-2}$ & $4.2\cdot 10^{-3}$
\end{tabular}
\caption{Numerical errors of the sharp-interface model with respect  to the $L^{2}(Q_T;\R^2)$ norm and with respect to the $L^{\infty}(Q_T;\R^2)$ norm for deformations upon change of time-step size $\{\tau,\tau/2,\tau/4\}$, uniform refinement of the mesh a) $\{h,\nicefrac{h}{2},\nicefrac{h}{4},\nicefrac{h}{8}\}$ in \Cref{fig:different_inital_meshes}, and polynomial order $\{P_1,P_2,P_3\}$ of the FE-space for $\bchi_{0}$.}
\label{tab:num_error_SI}
\end{table}

\section{Convergence to sharp-interface model}
\label{sec:convergence_discussion}
\Cref{fig:mobility_variation_DI_SI} gives a visual representation of the main results of this work and shows the direct comparison of phase-field model and sharp-interface model at time $t=T$ for different mobilities $m=m_0\varepsilon^\alpha$ and {$0\leq \alpha \le \infty$}. In the following we first give a detailed discussion of direct observations in this representation for different mobilities and then we discuss the (experimental) error $\|\bchi_{0,\Delta}-\bchi_{\varepsilon,\Delta'}\|$ in some detail for different norms.

\begin{figure}[hb!]\label{fig:convergence}
    \centering
    a)\,\includegraphics[height=0.20\textwidth,trim=0cm 38cm 4cm 28cm,clip]{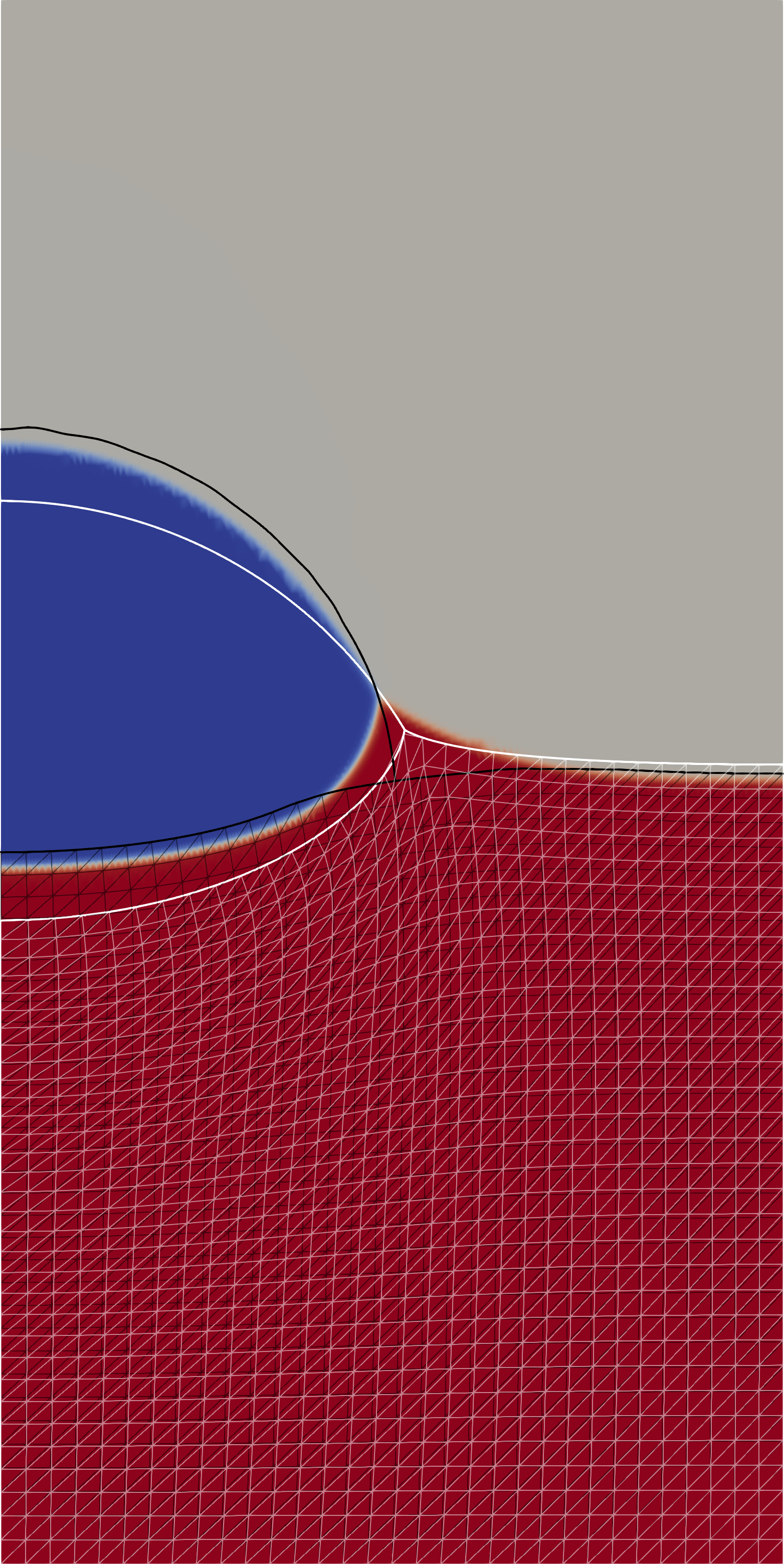}
    \includegraphics[height=0.20\textwidth,trim=18cm 50cm 22cm 42cm,clip]{pics/SI_DI_MInf.png}\hfill
    b)\,\includegraphics[height=0.20\textwidth,trim=0cm 38cm 4cm 28cm,clip]{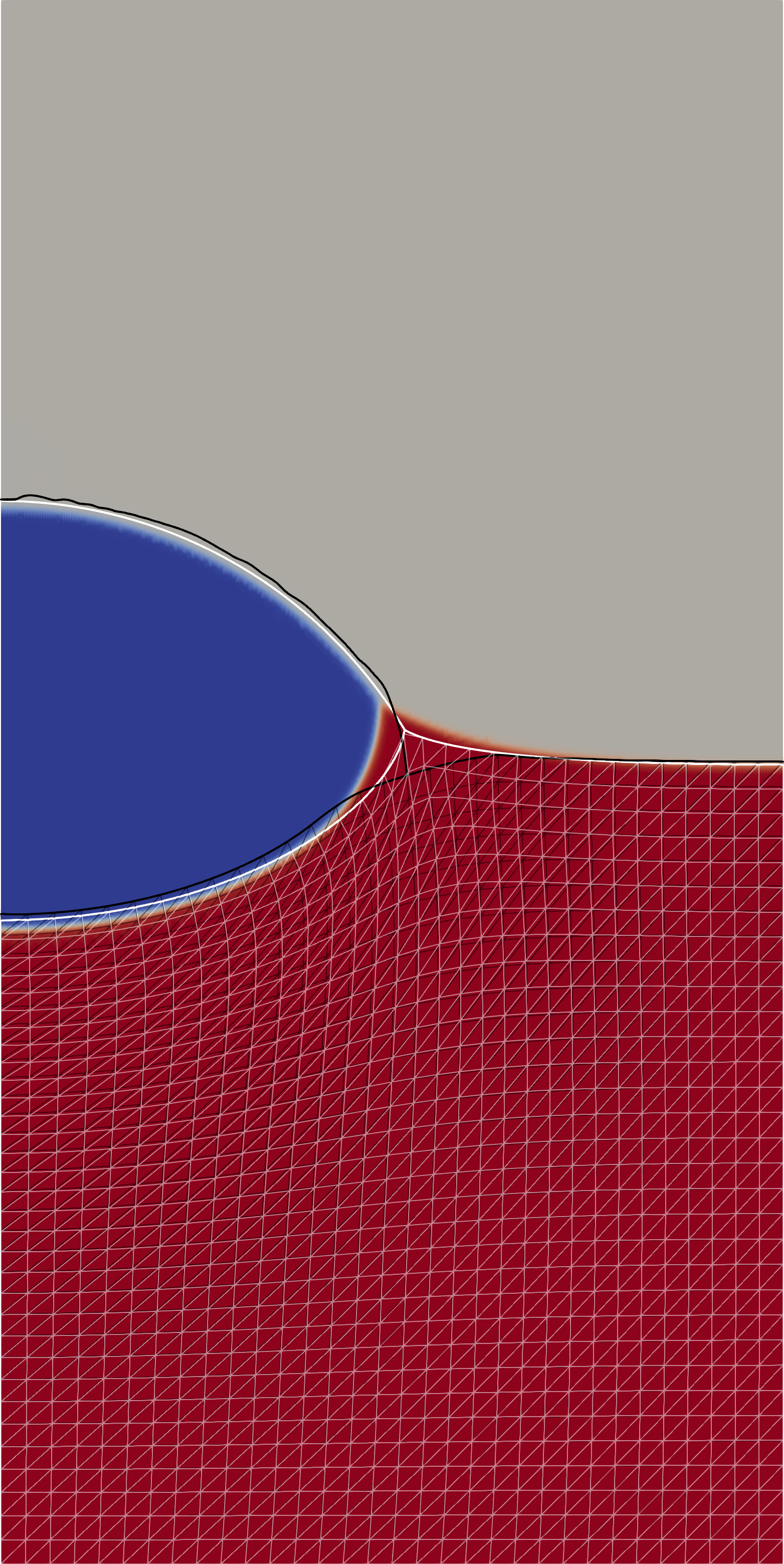}
    \includegraphics[height=0.20\textwidth,trim=18cm 50cm 22cm 42cm,clip]{pics/SI_DI_M1.png}\hfill
    
    c)\,\includegraphics[height=0.20\textwidth,trim=0cm 38cm 4cm 28cm,clip]{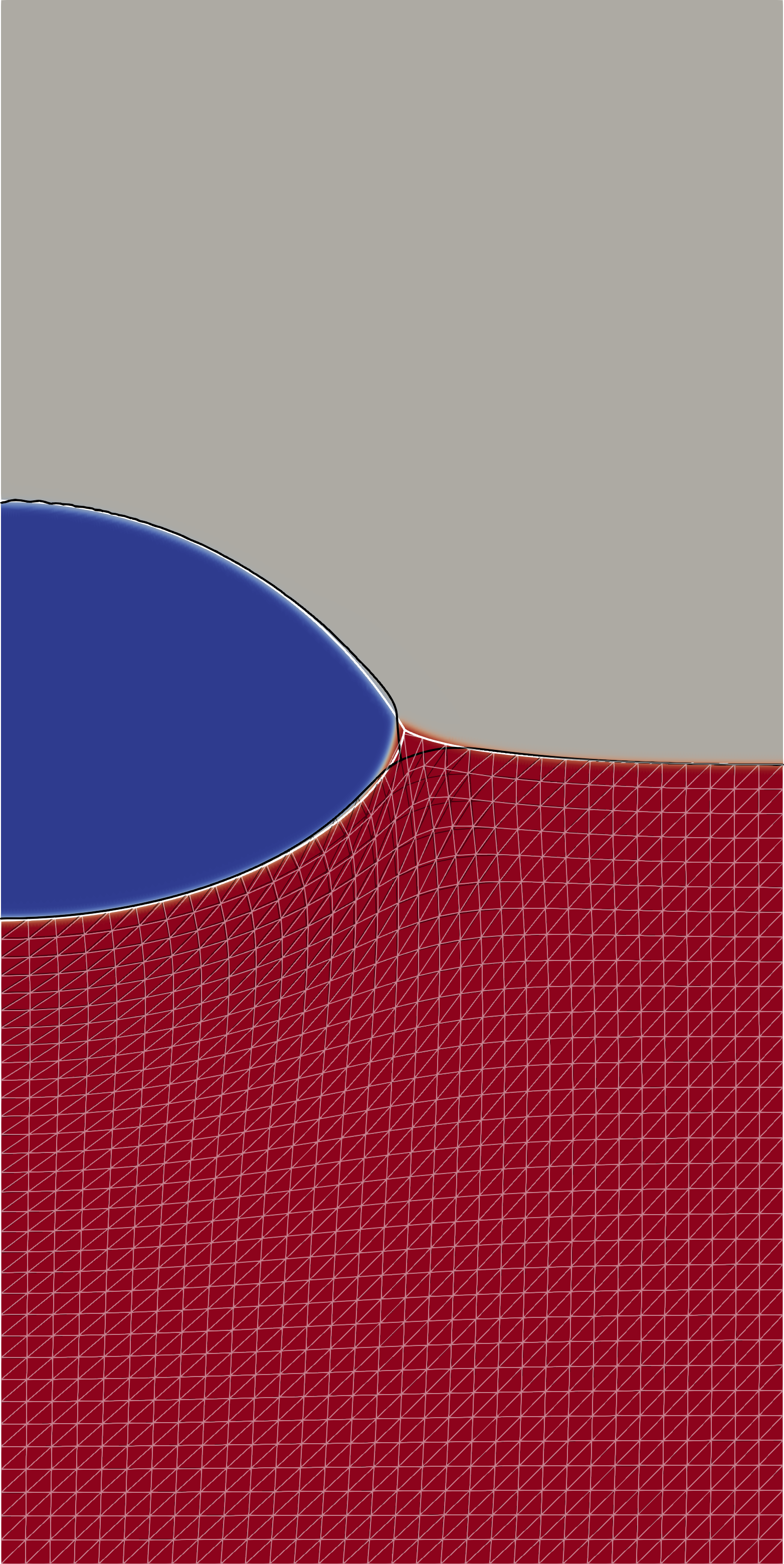}
    \includegraphics[height=0.20\textwidth,trim=18cm 50cm 22cm 42cm,clip]{pics/SI_DI_M15.png}\hfill
    d)\,\includegraphics[height=0.20\textwidth,trim=0cm 38cm 4cm 28cm,clip]{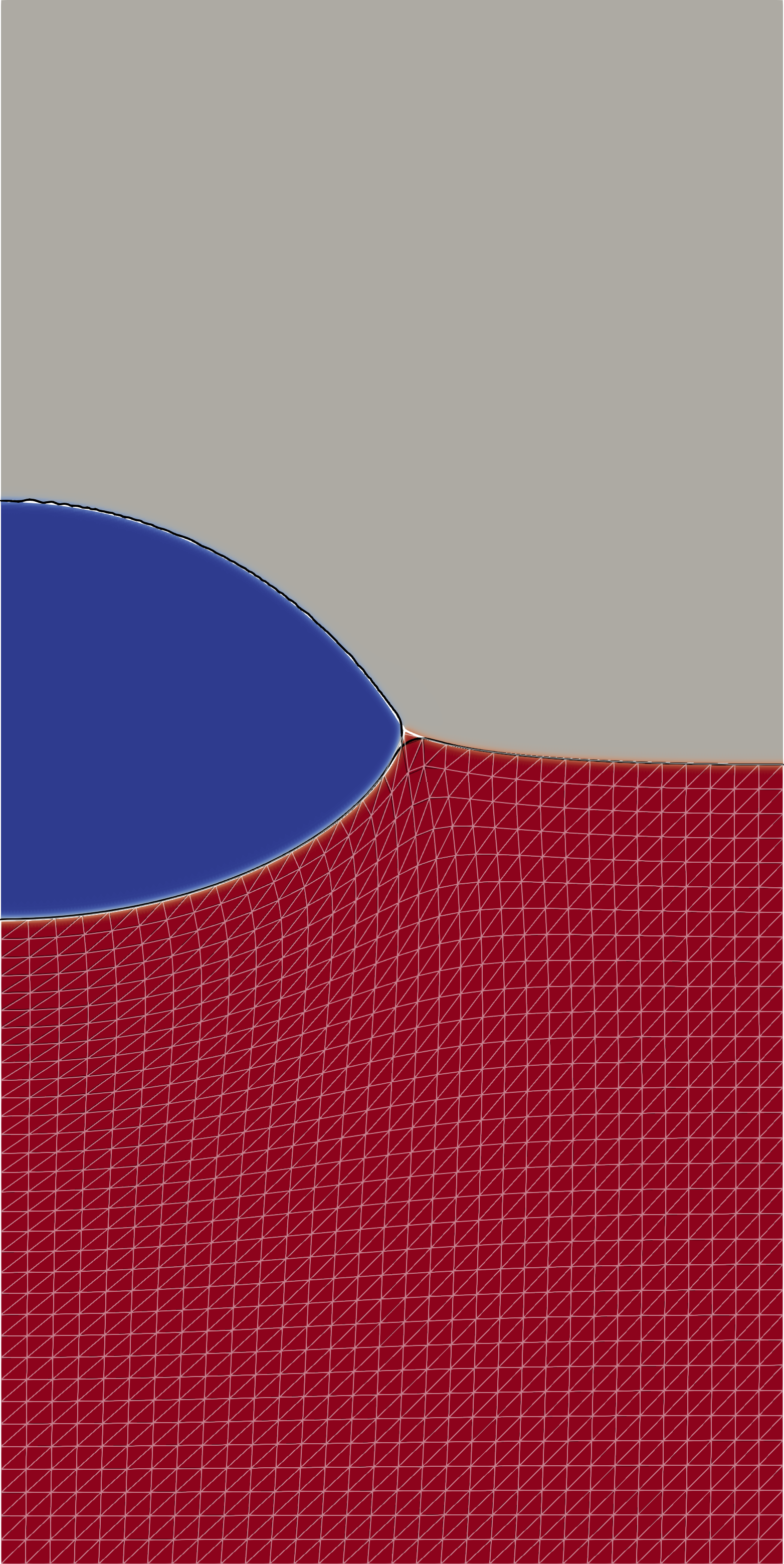}
    \includegraphics[height=0.20\textwidth,trim=18cm 50cm 22cm 42cm,clip]{pics/SI_DI_M2.png}
    
    e)\,\includegraphics[height=0.20\textwidth,trim=0cm 38cm 4cm 28cm,clip]{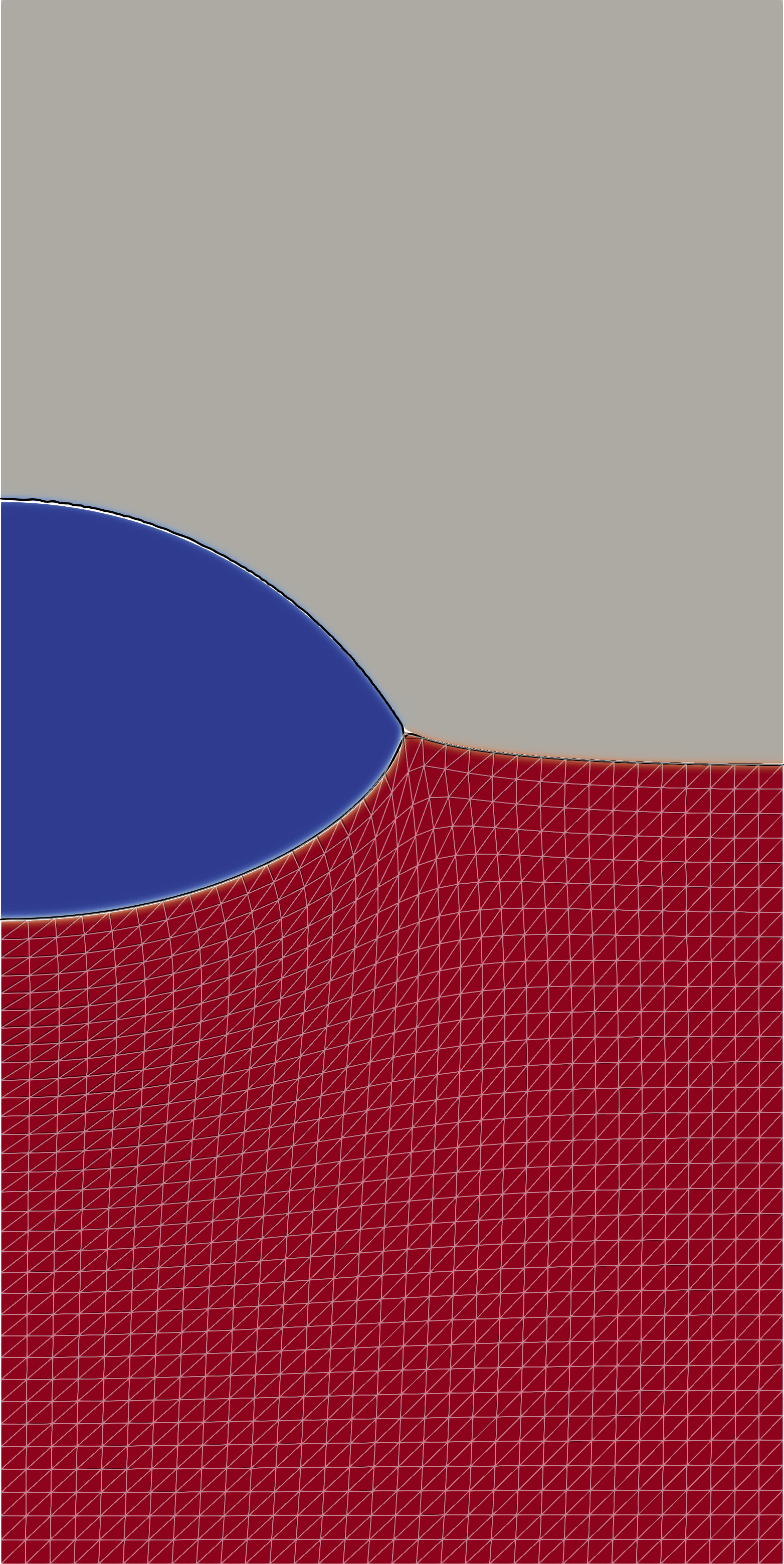}
    \includegraphics[height=0.20\textwidth,trim=18cm 50cm 22cm 42cm,clip]{pics/SI_DI_M25.png}\hfill
    f)\,\includegraphics[height=0.20\textwidth,trim=0cm 38cm 4cm 28cm,clip]{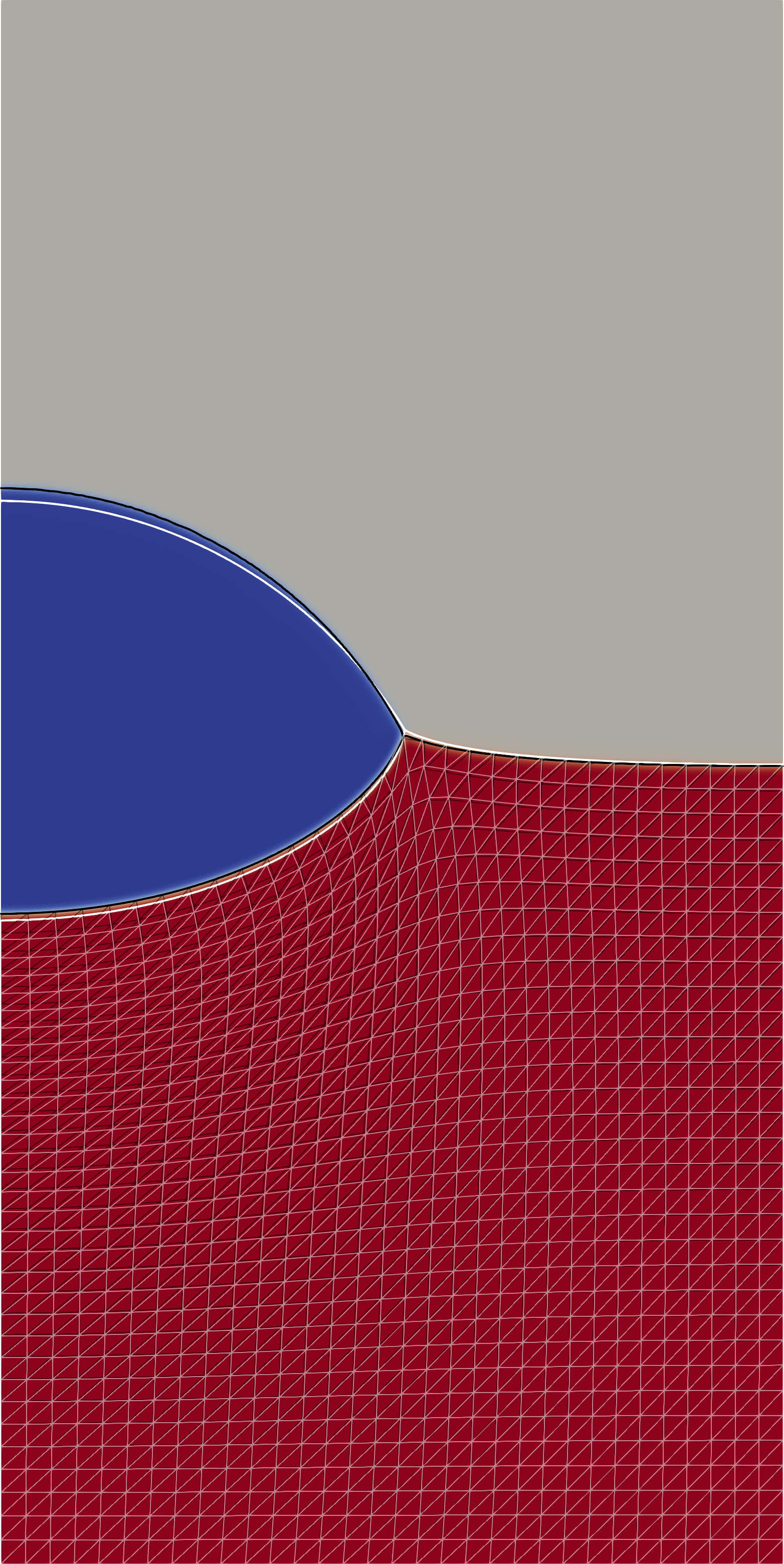}
    \includegraphics[height=0.20\textwidth,trim=18cm 50cm 22cm 42cm,clip]{pics/SI_DI_M3.png}
    
    g)\,\includegraphics[height=0.20\textwidth,trim=0cm 38cm 4cm 28cm,clip]{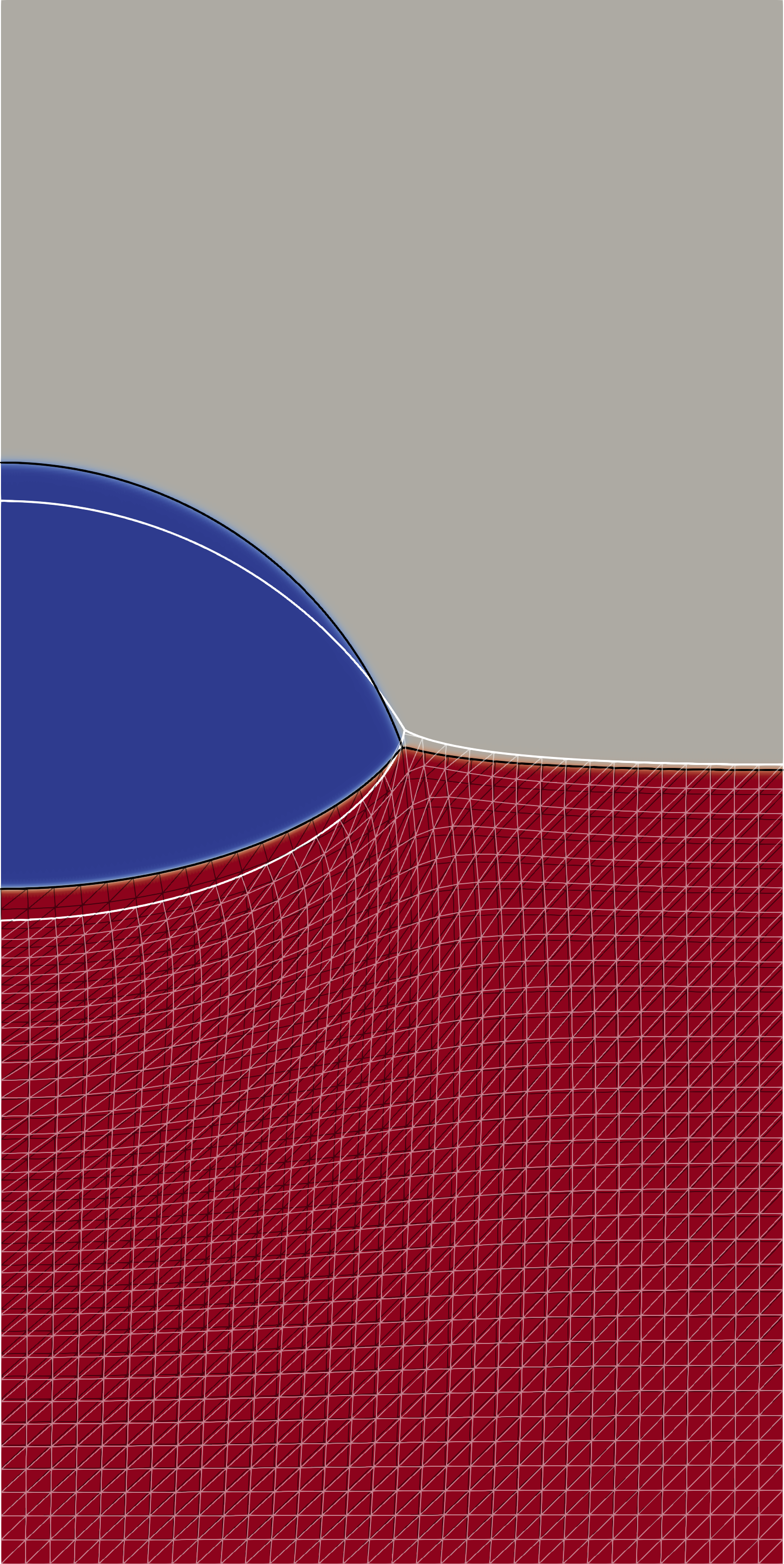}
    \includegraphics[height=0.20\textwidth,trim=18cm 50cm 22cm 42cm,clip]{pics/SI_DI_M4.png}\hfill
    h)\,\includegraphics[height=0.20\textwidth,trim=0cm 38cm 4cm 28cm,clip]{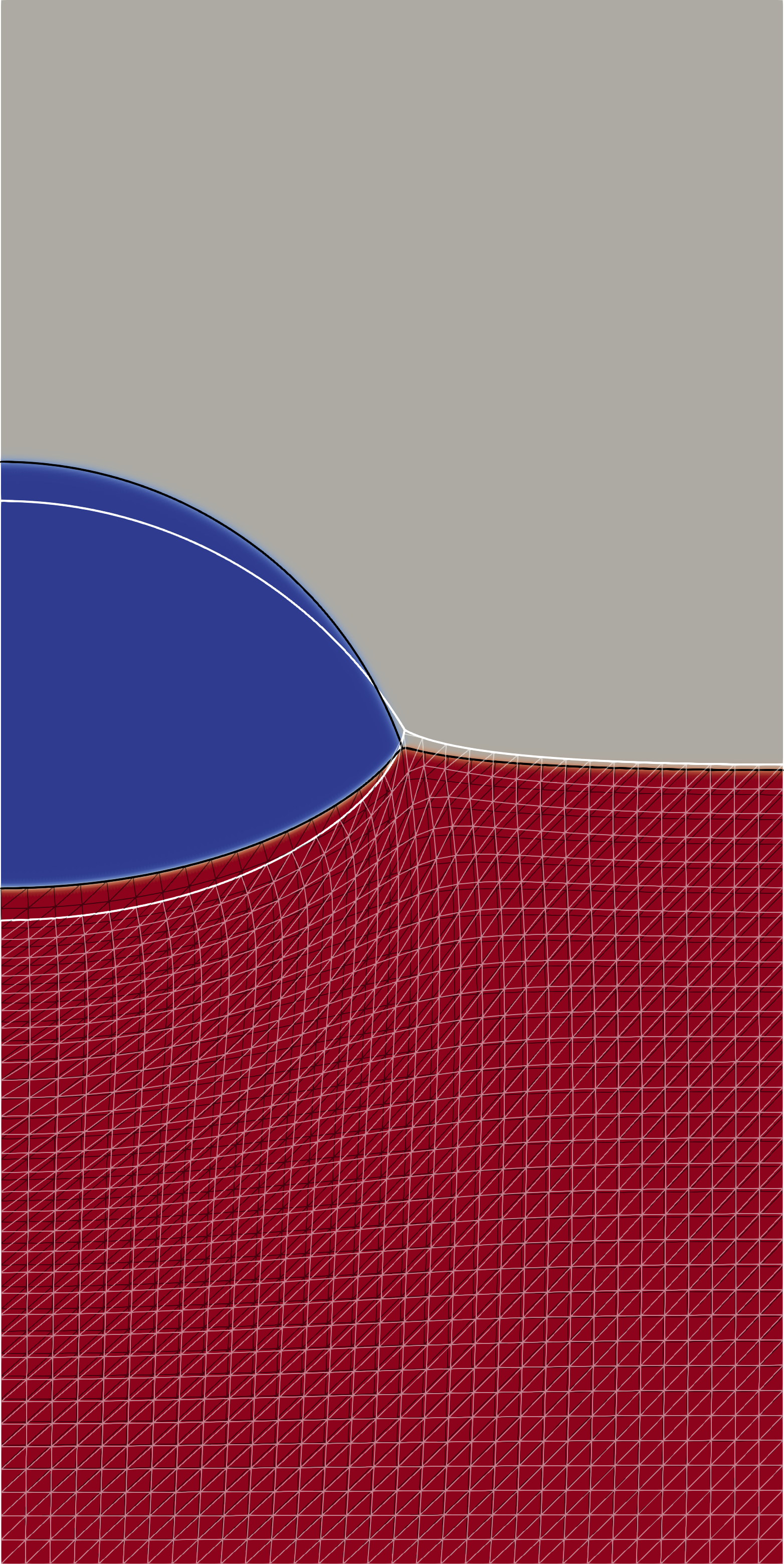}
    \includegraphics[height=0.20\textwidth,trim=18cm 50cm 22cm 42cm,clip]{pics/SI_DI_M0.png}

    \caption{Superposition of diffuse and sharp interfaces for $\varepsilon=2\bar{\varepsilon}$ with $k_{\bchi}=2$, $k_\psi=2$ and 2 uniform refinements of the mesh a) in \Cref{fig:initial_meshes}. The shading shows the phase-field indicator $\psi_{\rm id}$ with solid (red), liquid (blue) and air (gray) phases. Thick black lines display the phase-field deformation applied to the initial position of the diffuse interfaces and the thin black mesh displays the deformation of the initial diffuse domain of the elastic solid. Correspondingly, the thick white lines show the sharp interfaces $\bar{\Gamma}_{ij}(t)$ and the thin white mesh shows the displacement applied to the (sharp) elastic solid. We showcase different mobilities:     
    a) $\mob=1$ ($\alpha=0$),
    b) $\mob=\varepsilon$, c) $\mob=\varepsilon^{3/2}$, d) $\mob=\varepsilon^2$, e) $\mob=\varepsilon^{5/2}$, f) $\mob=\varepsilon^3$, g)  $\mob=\varepsilon^4$, h) $\mob=0$ ($\alpha=\infty$). Panels are shown at time $t=T$, except for a) which is at $t=0.1$. Right next to each image  a magnification of the contact line is shown.} 
    \label{fig:mobility_variation_DI_SI}
\end{figure}

Firstly, note that for mobility $m=1$, i.e. $\alpha=0$, in panel a) of \Cref{fig:mobility_variation_DI_SI} one can clearly see the nonconvergence of the phase-field model, since already at an earlier time the phase field (shading) and the displacement of the reference solid do not match. This is consistent with large excess diffusion leading to  a state where $\bF=\boldsymbol{I}$ and $\psi$ minimizing the phase-field energy separately, i.e. the evolution of the phase field disconnects from the evolution of the displacement and assumption A2 for the convergence is violated. This already happens early at time $t=0.1$.

Next, for mobilities $m=\varepsilon^\alpha$ for $\alpha\ge 1$ in panels {b)-e)} of \Cref{fig:mobility_variation_DI_SI} we see, except for slight convergent deviations in b), the phase field (shading) and the deformed initial interfaces (thick black lines) remain aligned, i.e. in those regions A2 is satisfied. However, for too small mobility, around $\alpha\ge 4$, the sharp interface and the diffuse field interface (thick white and black lines) may not be aligned as this additionally requires $\|\bchi_0-\bchi_\varepsilon\|\to 0$.

Further, note that even for $\alpha\ge 1$ in the cases b) $\alpha=1$, c) $\alpha=3/2$ we observe clear differences of the diffuse interface (shading) and the deformed mesh (thick black lines) much larger then the interfacial width $\varepsilon$ near the contact lines visible in the magnification of the contact line in \Cref{fig:mobility_variation_DI_SI}. What might even be counter-intuitive first is that the alignment of the sharp interface (thick white lines) and the diffuse interface (shading) in these cases b, c) appears much better then the alignment of the displacements (black and white thick lines). This again emphasizes that the convergence requires both convergence of the indicators A2 and convergence of the displacements A1. The main observation here is that while in a suitable (weak) norm one might prove convergence to the sharp-interface model, for all practical questions and in the strong $L^\infty$ norm, the error of the contact line position clearly lags behind the highly resolved interfacial thickness $\varepsilon$.

For the cases d) and e) with $\alpha=2$ and $\alpha=5/2$ we generally observe a very good agreement of phase fields (shading) and deformed initial diffuse interfaces (thick black lines) and deformed sharp interfaces (thick white lines) and the overall displacement (thin black and white mesh). This shows that these values are suitable for simulations when in particular a good precision of the predicted contact line position and small error in the displacements are desired.

For the cases f,g,h) with $\alpha=3,4,\infty$ we still see the perfect alignment of phase-field (shading) and deformed diffuse interfaces (thick black lines) but no convergence of the displacement, i.e. assumption A1 for convergence is violated. As explained earlier, this is due to the fact that the diffuse-interface profile does not reach the optimal $\tanh$-profile and the approximation of the Eulerian surface energy deteriorates.

{In our Lagrangian formulation}, it is possible to choose extremely small Cahn-Hilliard mobilities and even $m=0$ without the need to revise the numerical procedures and introduce artificial stabilization terms. While in an  Eulerian framework such a stabilization would be needed and possible restrictions due to CFL conditions might appear, the alignment of phase fields and displacement could in principle be also verified using an Eulerian approach by integrating the velocity field to a obtain a deformation map. This concludes the discussion of \Cref{fig:mobility_variation_DI_SI} and clearly shows that $\alpha=1$ produces significant errors near the contact line.

\begin{figure}
    \centering
   
    \includegraphics[width=0.49\textwidth]{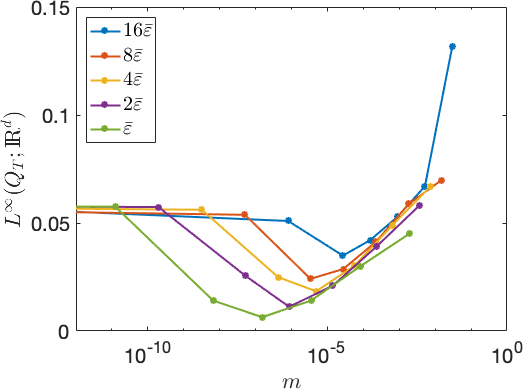}
    \includegraphics[width=0.49\textwidth]{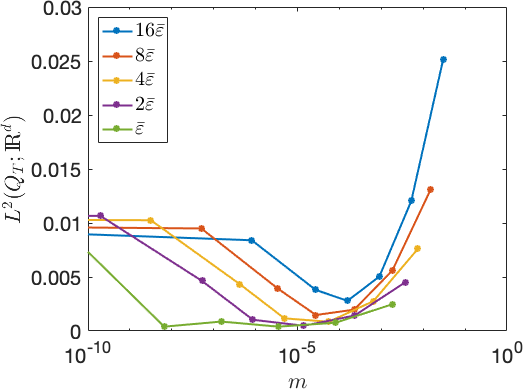}
    
    \caption{Convergence error $\|\bchi_{0,\Delta}-\bchi_{\varepsilon,\Delta'}\|$ for different $\varepsilon=2^n\bar\varepsilon$ and different mobilities $m(\varepsilon)=\varepsilon^\alpha$ as a function $m(\varepsilon)$. The used mobilities are $m=\varepsilon^{\infty}=0$, $m=\varepsilon^4$, $m=\varepsilon^3$, $m=\varepsilon^{5/2}$, $m=\varepsilon^2$, $m=\varepsilon^{3/2}$, $m=\varepsilon^1$, $m=1$. The left panel shows the convergence in the $L^\infty$ norm and the right panel in the $L^2$ Bochner norm.}
    \label{fig:norms_vs_mobility}
\end{figure}

The behavior shown in the spatial representation of the sharp-interface limit is mirrored in the behavior of the two norms shown in \Cref{fig:norms_vs_mobility}. Here we haven chosen a representation, where we plot different errors for different $\varepsilon=2^n\bar{\varepsilon}$ as a function of the mobility $m=\varepsilon^\alpha$. The reasons is that for any given $\varepsilon$ one would like to know the best possible mobility to obtain a sufficiently precise approximation of the sharp-interface model. From the preceding discussion its clear that both too large and too small values of the mobility will lead to suboptimal results. Both norms in \Cref{fig:norms_vs_mobility} clearly show that there are optimal values for $m(\varepsilon)$ for each $\varepsilon$ and this optimal mobility decreases for smaller $\varepsilon$. For $L^\infty$, while for $\varepsilon=16\bar\varepsilon$ we have $m_{\rm opt}\sim 10^{-4}$ we have for $\varepsilon=\bar\varepsilon$ an optimal value $\smash{m_{\rm opt}\sim 10^{-7}}$. This indicates that $\smash{m_{\rm opt}\sim \varepsilon^2}$ is optimal for a strong norm that takes errors near the contact line into consideration. On the other hand, the minimum in the $L^2$ error is rather broad and suggests that here, if the convergence in the contact line area is not of interest, a broader range of $\alpha$-values is possible.
Note that based on the discussion of numerical errors in 
\Cref{sec:num_errors}, these statements are limited by the numerical errors of order $10^{-3}$ in the $L^\infty$ norm. 

We finalize our discussion by considering the evolution of the spatial norms as a function of time $t$ on the cylinder $\Omega$ in \Cref{fig:errors_over_time}. For initial times $0<t<0.2$ we observe a steep increase in the convergence error, which might be affected due to the use of initial data that do not satisfy the Neumann triangle construction. This should result in a transient boundary layer at $t=0$ an possibly also influence the convergence rate of the sharp-interface limit in time. However, in our experience the relaxation of a droplet from a flat planar substrate is a more realistic scenario compared to well-prepared initial data that already satisfy the Neumann triangle. In the left panel of \Cref{fig:errors_over_time} we observe convergence in the $L^\infty$ norm for decreasing $\varepsilon=2^n\bar\varepsilon$. In the right panel we observe that for fixed $\varepsilon$ and mobility exponents $m=\varepsilon^\alpha$ with $\alpha=2,5/2,3$ we obtain the lowest errors of the order $10^{-2}$, which is still larger than the discretization errors previously discussed in \Cref{sec:num_errors}.

\begin{figure}
    \centering
   
    \includegraphics[width = 0.49\textwidth]{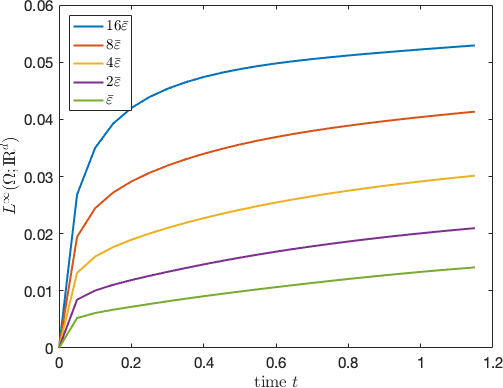}
    \includegraphics[width = 0.49\textwidth]{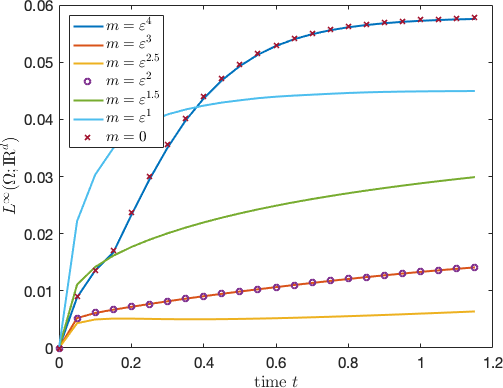}
    
    \caption{Convergence error $\|\bchi_{0,\Delta}-\bchi_{\varepsilon,\Delta'}\|(t)$ for $L^\infty(\Omega;\R^d)$ Sobolev
    norm as a function of time $t$. (Left) for fixed  mobility exponent $m=\varepsilon^2$ but different $\varepsilon$ and (right) for fixed $\varepsilon=\bar{\varepsilon}$ but different mobility exponents. }
    \label{fig:errors_over_time}
\end{figure}

\section{Conclusion}
We have presented a Lagrangian variational formulation of phase-field and sharp-interface models for a nonlinearly coupled fluid-structure interaction problem that involves nonlinear elasticity in a solid phase and a fluid and air phase, moving capillary interfaces and a moving contact line. We have discussed the sharp-interface limit depending on the mobility $m(\varepsilon)=m_0\varepsilon^\alpha$ in the Cahn-Hilliard model in different norms and shown that in the strong $L^\infty$ norm the moving contact lines suggests $\alpha_\mathrm{opt}=2$.
This discussion is also based on a detailed consideration of different contributions to numerical errors. Similar observations have been made before \cite{yue2010sharp,magaletti2013sharp}, but here we see them as a clear consequence of the presence of the moving contact line.

In the future one needs to avoid possible limitation due to locking by choosing incompressible models with suitable inf-sup stable finite elements. Additionally, by choosing degenerate state-dependent mobilities $m=m(\varepsilon,\psi)$ { and techniques as in \cite{dziwnik2017anisotropic}} one should be able to increase the region of validity compared to $m=m(\varepsilon)$. This discussion should be extended to models that take into account Navier-slip and dynamic contact angles. Rigorous and formal asymptotic results in this direction would be desirable but require a good understanding of both upper and lower bound $\underline{\alpha},\overline{\alpha}$. 